\documentclass[draft]{agujournal2019}
\usepackage{url} 
\usepackage{lineno}
\usepackage{bm}
\usepackage{upgreek}
\usepackage{amsmath}
\usepackage{amsfonts}
\usepackage{algorithm}
\usepackage{algorithmic}
\usepackage{color}
\usepackage{ragged2e}
\newcommand{\norm}[1]{\left\lVert#1\right\rVert}
\justifying

\draftfalse

\journalname{Water Resources Research}

\begin{document}

%
%


\title{Improving Simulation Efficiency of MCMC for Inverse Modeling of Hydrologic Systems with a Kalman-Inspired Proposal Distribution}

%
%




\authors{Jiangjiang Zhang\affil{1}, Jasper A. Vrugt\affil{2, 3}, Xiaoqing Shi\affil{4}, Guang Lin\affil{5}, Laosheng Wu\affil{6}, and Lingzao Zeng\affil{1}}


\affiliation{1}{Zhejiang Provincial Key Laboratory of Agricultural Resources and Environment, Institute of Soil and Water Resources and Environmental Science, College of Environmental and Resource Sciences, Zhejiang University, Hangzhou, China,}
\affiliation{2}{Department of Civil and Environmental Engineering, University of California, Irvine, California, USA,}
\affiliation{3}{Department of Earth System Science, University of California, Irvine, California, USA,}
\affiliation{4}{Key Laboratory of Surficial Geochemistry of Ministry of Education, School of Earth Sciences and Engineering, Nanjing University, Nanjing, China,}
\affiliation{5}{Department of Mathematics and School of Mechanical Engineering, Purdue University, West Lafayette, Indiana, USA,}
\affiliation{6}{Department of Environmental Sciences, University of California, Riverside, California, USA.}





\correspondingauthor{J. A. Vrugt}{jasper@uci.edu}
\correspondingauthor{L. Zeng}{lingzao@zju.edu.cn}




\begin{keypoints}
\item MCMC methods remain rather inefficient in exploring high-dimensional target distributions
\item This paper introduces a Kalman-inspired proposal distribution to improve the simulation efficiency of MCMC
\item The Kalman-inspired proposal distribution can be conveniently embedded in any adequate MCMC method
\end{keypoints}

%
%


\begin{abstract}
Bayesian analysis is widely used in science and engineering for real-time forecasting, decision making, and to help unravel the processes that explain the observed data. These data are some deterministic and/or stochastic transformations of the underlying parameters. A key task is then to summarize the posterior distribution of these parameters. When models become too difficult to analyze analytically, Monte Carlo methods can be used to approximate the target distribution. Of these, Markov chain Monte Carlo (MCMC) methods are particularly powerful. Such methods generate a random walk through the parameter space and, under strict conditions of reversibility and ergodicity, will successively visit solutions with frequency proportional to the underlying target density. This requires a proposal distribution that generates candidate solutions starting from an arbitrary initial state. The speed of the sampled chains converging to the target distribution deteriorates rapidly, however, with increasing parameter dimensionality. In this paper, we introduce a new proposal distribution that enhances significantly the efficiency of MCMC simulation for highly parameterized models. This proposal distribution exploits the cross-covariance of model parameters, measurements and model outputs, and generates candidate states much alike the analysis step in the Kalman filter. We embed the Kalman-inspired proposal distribution in the DREAM algorithm during burn-in, and present several numerical experiments with complex, high-dimensional or multi-modal target distributions. Results demonstrate that this new proposal distribution can greatly improve simulation efficiency of MCMC. Specifically, we observe a speed-up on the order of 10 - 30 times for groundwater models with more than one-hundred parameters.
\end{abstract}

%
%

%


%
%
%
%

\section{Introduction and Scope}
Mathematical modeling has become an integral part of the scientific method. Computer simulation is particularly useful for hypothesis testing, decision making, to gain (new) insights and understanding of system functioning, and to predict system behavior into the space (interpolation) and time (forecasting) domain. The complexity of hydrologic systems poses significant modeling challenges, in particular how to characterize adequately water flow and storage in the presence of (often) incomplete and insufficient observations, process knowledge and system characterization. This necessitates a systematic framework for uncertainty quantification of model simulations. These uncertainties originate from model conceptualization and data collection, and include model structural errors, measurement errors of the initial conditions, forcing data and model output data, and uncertainty in the model parameters \cite{kavetski2006a,kavetski2006b,refsgaard2012,renard2011,vrugt2005,wagener2005,zheng2018}. 

Let us assume that the data-generating process of some arbitrary real-world system can be written as
\begin{linenomath*}
\begin{equation}
\widetilde{\textbf{d}} = f(\bm{\uptheta}) + \bm{\upepsilon},
\label{model}
\end{equation}
\end{linenomath*}
where $\widetilde{\textbf{d}} = \{ \widetilde{d}_1,\ldots,\widetilde{d}_n \}$ is a $n$-vector of measurements, $f(\cdot)$ signifies a computer model of the system of interest, $\bm{\uptheta} = \{\theta_{1},\ldots,\theta_{k} \}$ denotes a $k$-vector of model parameters, and $\bm{\upepsilon} = \{ \varepsilon_1,\ldots,\varepsilon_{n} \}$ represents a $n$-vector of measurement errors, respectively. The Bayesian paradigm treats the parameters in equation (\ref{model}) as random variables with joint probability density function. This multivariate distribution, the so-called posterior parameter distribution, $p(\bm{\uptheta}|\widetilde{\textbf{d}})$, is the consequence of two antecedents, a prior distribution, $p(\bm{\uptheta})$, which captures our initial degree of beliefs in the values of the model parameters, and a likelihood function, $L(\bm{\uptheta}|\widetilde{\textbf{d}}) \equiv p(\widetilde{\textbf{d}} | \bm{\uptheta})$, which quantifies by the rules of probability theory the level of confidence in the parameter values, $\bm{\uptheta}$, in light of the observed data, $\widetilde{\textbf{d}}$. Bayes' theorem expresses mathematically the relationship between the prior, conditional and posterior distribution of the parameters, $\bm{\uptheta}$, as follows
\begin{linenomath*}
\begin{equation}
p(\bm{\uptheta} | \widetilde{\textbf{d}}) = \frac{p(\bm{\uptheta}) p(\widetilde{\textbf{d}} | \bm{\uptheta})}{ \int p(\bm{\uptheta}) p(\widetilde{\textbf{d}} | \bm{\uptheta}) \text{d} \bm{\uptheta} } \propto p(\bm{\uptheta}) L(\bm{\uptheta}|\widetilde{\textbf{d}}),
\label{bayes}
\end{equation}
\end{linenomath*}
where the denominator, $p(\widetilde{\textbf{d}}) = \int p(\bm{\uptheta}) p(\widetilde{\textbf{d}} | \bm{\uptheta}) \text{d} \bm{\uptheta}$, is the so-called Bayesian evidence, evidence, marginal likelihood or model likelihood. This normalization constant guarantees that $p(\bm{\uptheta} | \widetilde{\textbf{d}})$ integrates to unity, i.e., $\int_{\bm{\Theta}} p(\bm{\uptheta} | \widetilde{\textbf{d}}) \text{d} \bm{\uptheta} = 1$, where $\bm{\Theta}\subseteq \mathbb{R}^{k}$ signifies the feasible parameter space, and $p(\bm{\uptheta} | \widetilde{\textbf{d}}) \geq 0 \mbox{ for all } \bm{\uptheta} \in \bm{\Theta}$. The evidence is of great importance for hypothesis testing via model selection \cite{cao2018,volpi2017,zeng2018}, but can be discarded if our interest lies in estimation of the posterior distribution of the parameters. Henceforth, one often removes the denominator from Bayes' theorem and works instead with the unnormalized density, i.e., the right-hand-side of equation (\ref{bayes}). When the measurement errors, $\bm{\upepsilon}$, are normally distributed with zero-mean and $n \times n$ covariance matrix, $\textbf{R}$, i.e., $\bm{\upepsilon} \sim \mathcal{N}_n(\textbf{0},\textbf{R})$, the likelihood function can be expressed as
\begin{linenomath*}
\begin{equation}
L(\bm{\uptheta}|\widetilde{\textbf{d}}) = \frac{1}{(2\pi)^{n/2}\left| \textbf{R} \right|^{1/2}} \exp \left\{ -\frac{1}{2} \left [ \widetilde{\textbf{d}} - f(\bm{\uptheta}) \right ]^\text{T} \textbf{R}^{-1} \left[ \widetilde{\textbf{d}} - f(\bm{\uptheta}) \right] \right\},
\label{likelihood}
\end{equation}
\end{linenomath*}
where $\left| \cdot \right|$ signifies the determinant operator. Nevertheless, in some cases, one has to consider correlated, heteroscedastic and non‐Gaussian measurement errors. Then other forms of likelihood functions can be adopted, e.g., the formal likelihood function proposed by \citeA{schoups2010}.

A key task is now to summarize the posterior parameter distribution, $p(\bm{\uptheta}|\widetilde{\textbf{d}})$. In most practical cases, $p(\bm{\uptheta} | \widetilde{\textbf{d}})$ cannot be derived by analytical means nor by analytical approximation, and Monte Carlo methods can be used to sample the target distribution. Of these, Markov chain Monte Carlo (MCMC) methods have become increasingly popular in the past decades. Such methods generate a (quasi-)random walk through the parameter space and collect the sampled solutions in one or more Markov chains. Under strict conditions of reversibility and ergodicity, the chain(s) will gradually converge to an equilibrium distribution equivalent to the target distribution. This means that if one looks at the archived values of $\bm{\uptheta}$ in the chain(s) sufficiently far from the arbitrary initial state(s), i.e., after the so-called burn-in period, then these successively generated states will be distributed according to $p(\bm{\uptheta}|\widetilde{\textbf{d}})$, the unknown target distribution of $\bm{\uptheta}$. In practice, it is necessary to monitor convergence of the sampled chain(s) with some diagnostic metrics \cite{brooks1998}. Then we discard the Markov chain states in the burn-in period to minimize the effect of initial values on the posterior inference, and use the remaining states to derive the desired statistics of $p(\bm{\uptheta}|\widetilde{\textbf{d}})$.

The earliest MCMC method, i.e., the random walk Metropolis (RWM) algorithm, was developed by \citeA{metropolis1953} and works as follows. First, a candidate, $\bm{\uptheta}_\text{p}$, is sampled from a symmetric proposal distribution, $q(\cdot)$, which is centered around the current state, $\bm{\uptheta}_{(t-1)}$, of the Markov chain. Next, the candidate is accepted with Metropolis probability
\begin{linenomath*}
\begin{equation}
p_{\text{acc}}(\bm{\uptheta}_{(t-1)} \to \bm{\uptheta}_\text{p}) = \min \left [1,\frac {p(\bm{\uptheta}_\text{p}) L(\bm{\uptheta}_\text{p} | \widetilde{\textbf{d}})} { p(\bm{\uptheta}_{(t-1)}) L(\bm{\uptheta}_{(t-1)} | \widetilde{\textbf{d}})} \right] = \min \left[1,\frac {p(\bm{\uptheta}_\text{p} | \widetilde{\textbf{d}})} {p(\bm{\uptheta}_{(t-1)} | \widetilde{\textbf{d}})} \right].
\label{RWM}
\end{equation}
\end{linenomath*}
Finally, if the candidate is accepted, the chain moves to $\bm{\uptheta}_\text{p}$, otherwise the chain remains at its current state, $\bm{\uptheta}_{(t-1)}$. Repeated application of these three steps results in a Markov chain with equilibrium distribution, $p(\bm{\uptheta}|\widetilde{\textbf{d}})$. \citeA{hastings1970} generalized the RWM algorithm to asymmetric proposal distributions when the probability density of the forward jump, $q(\bm{\uptheta}_\text{p}|\bm{\uptheta}_{(t-1)})$, does not equal the probability density of the backward jump, $q(\bm{\uptheta}_{(t-1)}|\bm{\uptheta}_\text{p})$. A simple correction of the acceptance probability will neutralize this imbalance in jump probabilities
\begin{linenomath*}
\begin{equation}
p_{\text{acc}}(\bm{\uptheta}_{(t-1)} \to \bm{\uptheta}_\text{p}) = \min \left[1,\frac{p(\bm{\uptheta}_\text{p} |\widetilde{\textbf{d}})}{p(\bm{\uptheta}_{(t-1)} | \widetilde{\textbf{d}} )} \cdot \frac{q(\bm{\uptheta}_{(t-1)}|\bm{\uptheta}_\text{p})} { q(\bm{\uptheta}_\text{p}|\bm{\uptheta}_{(t-1)})} \right],
\label{MH}
\end{equation}
\end{linenomath*}
where $q(\bm{\uptheta}_{(t-1)}|\bm{\uptheta}_\text{p})$ and $q(\bm{\uptheta}_\text{p}|\bm{\uptheta}_{(t-1)})$ signify the conditional probabilities of trail moves from $\bm{\uptheta}_\text{p}$ to $\bm{\uptheta}_{(t-1)}$ and from $\bm{\uptheta}_{(t-1)}$ to $\bm{\uptheta}_\text{p}$, respectively. Equation (\ref{MH}) is also known as the Metropolis-Hastings (MH) algorithm and has enabled the practical application of Bayesian inference to a very large class of models and data. Indeed, the MH algorithm has revolutionized the field of computational statistics, and because of this reason was elected as one of the top 10 most important algorithms of the 20th century \cite{beichl2000}.

The efficiency of the RWM and MH algorithms depends in large part on the scale and orientation of the proposal distribution, $q(\cdot)$, used to create trial moves (transitions) in the Markov chain. When the proposal distribution is too dispersed, a disproportionately large number of candidate states will be rejected, and the chain will only converge slowly to the target distribution. On the other hand, when the proposal distribution is too narrow (under-dispersed), most candidate states will be accepted, but the chain will also travel slowly to the target distribution as the update is very small, and MCMC will experience a very long burn-in period.

Much research has been devoted in the past decades to improving the efficiency of MCMC methods \cite{brooks2011,calderhead2014,gilks1996}. This includes the use of kernel adaptation \cite{gilks1994,haario1999,haario2001,kuczera1998,vrugt2003}, Hamiltonian dynamics \cite{duane1987,hoffman2014,neal2011}, multi-try proposals \cite{laloy2012,liu2000}, Langevin dynamics \cite{girolami2011,roberts1998,roberts2002}, parallel marginalization \cite{weare2007}, delayed-rejection \cite{haario2006}, bacterial kernels \cite{yang2013}, kernel coupling and multiple chain simulation \cite{craiu2009,terbraak2006,terbraak2008,vrugt2008b,vrugt2009}, early rejection \cite{laloy2013,solonen2012}, parallelization \cite{calderhead2014,neiswanger2013} and pre-fetching \cite{brockwell2006,strid2010}. These algorithms should, at least in theory, converge asymptotically to the target distribution. Yet, their performance can differ tremendously over a finite number of model evaluations.

Notwithstanding these methodological advances, MCMC methods still remain rather inefficient in exploring high-dimensional parameter spaces. This severely impairs their practical application to complex system models. The goal of this paper is twofold. First, we introduce a new proposal distribution that enhances significantly the efficiency of MCMC simulation for highly parameterized system models. This proposal distribution was suggested in \citeA{vrugt2013} and generates candidate states much alike the analysis step in the Kalman filter \cite{kalman1960}. However, \protect\citeA{vrugt2013} only foresaw the potential of this new idea and it has never been actually implemented and systematically examined. Second, to turn this raw idea into a reliable product, we incorporate the Kalman-inspired proposal distribution with a selection probability, $p_{\rm{K}}\in(0,1)$, in various MCMC algorithms, especially the DiffeRential Evolution Adaptive Metropolis (DREAM) algorithm \cite{vrugt2016}. The DREAM algorithm has shown to work well in a large array of inverse problems, involving high-dimensional and/or multi-modal target distributions with complex multivariate parameter dependencies \cite{bikowski2012,laloy2012,laloy2013,muleta2012,ramin2014,shi2012,shi2014,wohling2011,zhang2018}. As the Kalman-inspired proposal can introduce asymmetry to the sampled chains, it is suggested to restrict its use only during burn-in (recommended), or use the MH algorithm or a randomized move direction to neutralize the asymmetry. The resulting DREAM$_\text{(KZS)}$ algorithm incorporates the Kalman-inspired proposal distribution and is a sibling of the DREAM family \cite{vrugt2016}.

The remainder of this paper is organized as follows. In section 2, we introduce elements and theory of the Kalman-inspired proposal distribution, and provide a detailed recipe and discussion of the DREAM$_\text{(KZS)}$ algorithm. Section 3 presents the results of the DREAM$_\text{(KZS)}$ algorithm for several numerical studies with complex, high-dimensional or multi-modal target distributions. In this section we are especially concerned with benchmark analysis of the DREAM$_\text{(KZS)}$ algorithm. Moreover, to demonstrate that the Kalman-inspired proposal distribution can be conveniently embedded in any adequate MCMC method, we further test its performance in two plainer MCMC algorithms \cite{haario2001,haario2006}. Finally, we conclude this paper and provide some further discussion in section 4.

\section{Methods}
In this section, we introduce the Kalman-inspired proposal distribution and discuss its implementation in the DREAM algorithm that adopts a mix of parallel direction and snooker proposal distributions, i.e., the DREAM$_\text{(ZS)}$ algorithm \cite{vrugt2016}.

\subsection{Markov Chain Monte Carlo Simulation}
The core of the MH algorithm can be written in just a few lines (see Algorithm \ref{alg:RWM}). This algorithm simulates a single chain trajectory using the proposal distribution, $q(\cdot)$, a random number generator on the $(0,1)$ interval, $z \sim \mathcal{U}(0,1)$, and the target density, $p(\bm{\uptheta}|\widetilde{\textbf{d}})$ in equation (\ref{bayes}).
\begin{algorithm}[H]
	\caption{Metropolis-Hastings}
	\begin{algorithmic}[1]
		\STATE Define the desired length of the Markov chain, $T$.
		\STATE Draw the initial chain state, $\bm{\uptheta}_{(1)}$, from the prior distribution, $\bm{\uptheta}_{(1)} \sim p(\bm{\uptheta})$.
		\STATE Evaluate the unnormalized target density, $p(\bm{\uptheta}_{(1)}|\widetilde{\textbf{d}}) = p(\bm{\uptheta}_{(1)})L(\bm{\uptheta}_{(1)}|\widetilde{\textbf{d}})$, at $\bm{\uptheta}_{(1)}$.
		\FOR{$t = 2,\ldots,T$}
		\STATE Sample a candidate, $\bm{\uptheta}_\text{p}$, from a proposal distribution, $q(\cdot|\bm{\uptheta}_{(t-1)})$, $\bm{\uptheta}_\text{p} \sim q(\bm{\uptheta}|\bm{\uptheta}_{(t-1)})$. 
		\STATE Compute the unnormalized target density, $p(\bm{\uptheta}_\text{p}|\widetilde{\textbf{d}}) = p(\bm{\uptheta}_\text{p}) L(\bm{\uptheta}_\text{p} | \widetilde{\textbf{d}})$, at $\bm{\uptheta}_\text{p}$.
		\STATE Calculate the acceptance probability, $p_{\text{acc}}(\bm{\uptheta}_{(t-1)} \to \bm{\uptheta}_\text{p})$, in equation (\ref{MH}).
		\STATE Draw a random sample, $z \sim \mathcal{U}(0,1)$.
		\IF{$p_{\text{acc}}(\bm{\uptheta}_{(t-1)} \to \bm{\uptheta}_\text{p})\ge z$}
		\STATE Set $\bm{\uptheta}_{(t)} = \bm{\uptheta}_\text{p}$,
		\ELSE
		\STATE Set $\bm{\uptheta}_{(t)} = \bm{\uptheta}_{(t-1)}$.
		\ENDIF
		\ENDFOR
	\end{algorithmic}
	\label{alg:RWM}
\end{algorithm}
The user is free to select the proposal distribution as the acceptance probability in equation (\ref{MH}) preserves the underlying stationary distribution of the Markov chain. A common choice is the multivariate normal distribution, $q = \mathcal{N}_{k}(\cdot,s_{k}\bm{\Sigma})$, with covariance matrix, $\bm{\Sigma}$, and scaling factor, $s_{k}$. To enhance search efficiency and protect against an inadequate selection of the proposal distribution, we can update the scale and orientation of the proposal distribution every $m \gg 1$ iterations using all past samples stored in the chain, $\bm{\Sigma} = {\rm Cov}\left(\{\bm{\uptheta}_{(1)},\ldots,\bm{\uptheta}_{(t-1)}\}\right) + \varphi \textbf{I}_{k}$, where $\textbf{I}_{k}$ denotes the $k \times k$ identity matrix, and $\varphi = 10^{-6}$ is a small scalar that prevents the collapse of the sample covariance matrix to singularity (jumps become zero). This is the so-called Adaptive Metropolis (AM) algorithm of \citeA{haario2001}. As a basic choice, the scaling factor is often chosen to be $s_{k} = 2.38^{2}/{k}$, which has proven optimal for Gaussian target and proposal distributions \cite{gelman1996,roberts1997}, and should result in an acceptance rate close to $0.44$ for $k = 1$, $0.28$ for $k = 5$ and $0.23$ for a large $k$. It should be evident that the choice of the proposal distribution determines in large part the sampling efficiency and convergence speed of the Markov chain. What would be desirable is a proposal distribution that can transit the chain to the target distribution quickly.

\subsection{The Kalman-Inspired Proposal Distribution}
This section introduces an alternative proposal distribution designed to accelerate the movement of the chain to the posterior distribution, $p(\bm{\uptheta}|\widetilde{\textbf{d}})$. In this section, we present the theory and formulas of the so-called Kalman-inspired proposal distribution.

This new proposal distribution is inspired by the state analysis step in the Kalman filter \cite{kalman1960}. The Kalman-inspired proposal distribution uses explicitly the distance of the $n$-vector of simulated model outputs, $\textbf{d} = f(\bm{\uptheta})$, to the observed data, $\widetilde{\textbf{d}}$,
\begin{linenomath*}
\begin{equation}
\begin{split}
\bm{\uptheta}_\text{p} & = \bm{\uptheta}_{(t-1)} + \textbf{C}_{\bm{\uptheta}\textbf{d}}( \textbf{C}_{\textbf{d}\textbf{d}} + \textbf{R} )^{-1} \left[ \widetilde{\textbf{d}} + \bm{\upepsilon}_{(t-1)}- f(\bm{\uptheta}_{(t-1)}) \right] \\
& = \bm{\uptheta}_{(t-1)} + \textbf{K} \textbf{r}_{(t-1)} + \textbf{K} \bm{\upepsilon}_{(t-1)}\\
& = \bm{\uptheta}_{(t-1)} + \Delta \bm{\uptheta}_{(t-1)},
\end{split}
\label{Kalman_jump}
\end{equation}
\end{linenomath*}
where $\textbf{C}_{\bm{\uptheta}\textbf{d}} = \text{Cov}(\bm{\uptheta},\textbf{d})$ denotes the $k \times n$ cross-covariance matrix of model parameters and model outputs, $\textbf{C}_{\textbf{d}\textbf{d}} = \text{Cov}(\textbf{d},\textbf{d})$ signifies the $n \times n$ auto-covariance matrix of model outputs, $\textbf{R}$ is the $n \times n$ covariance matrix of measurement errors in equation (\ref{likelihood}), $\bm{\upepsilon}_{(t-1)}$ is a $n$-vector of random draw from the distribution of measurement errors, $\bm{\upepsilon}_{(t-1)} \sim \mathcal{N}_{n}(\textbf{0},\textbf{R})$, $\textbf{r}_{(t-1)} = \widetilde{\textbf{d}} - f(\bm{\uptheta}_{(t-1)})$ is the residual vector of parameters $\bm{\uptheta}_{(t-1)}$, and $\textbf{K}$ constitutes the so-called Kalman gain
\begin{linenomath*}
\begin{equation}
\textbf{K} = \textbf{C}_{\bm{\uptheta}\textbf{d}}( \textbf{C}_{\textbf{d}\textbf{d}} + \textbf{R} )^{-1}.
\label{Kalman}
\end{equation}
\end{linenomath*}
It is noted here that, the Kalman-inspired proposal uses covariances that converge, after burn in, towards posterior covariances. Besides, this proposal is only used to generate a candidate point, $\bm{\uptheta}_\text{p}$, from the current chain state, $\bm{\uptheta}_{(t-1)}$, not to directly estimate the posterior distribution, $p(\bm{\uptheta} | \widetilde{\textbf{d}})$. In other words, the purpose and implementation of the Kalman-inspired proposal distribution are quite different from the Kalman filter, or its Monte Carlo variants.

The jump vector of the Kalman-inspired proposal, $\Delta \bm{\uptheta}_{(t-1)}$, in equation (\ref{Kalman_jump}) is made up of two different components: A deterministic displacement vector, $\textbf{K} \textbf{r}_{(t-1)}$, which orients the jump towards the ``true" parameters, $\bm{\uptheta}^{*}$, of the data-generating process, and a random displacement vector, $\textbf{K} \bm{\upepsilon}_{(t-1)}$, with zero-mean normally distributed variables, which introduce randomness into the sampled candidate states. This latter term of the jump vector, $\Delta \bm{\uptheta}_{(t-1)}$, is of crucial importance as it enables the Kalman-inspired proposal to sample, with non-zero probability, all possible states of the target distribution. This ensures ergodicity of the sampled chain. From equation (\ref{Kalman_jump}), it is clear that the width of the Kalman-inspired proposal distribution, which is determined by the random displacement part, $\textbf{K} \bm{\upepsilon}_{(t-1)}$, should keep relatively stable with time, as the distribution of measurement errors is constant. On the contrary, the deterministic jump distance, $\textbf{K} \textbf{r}_{(t-1)}$, should decrease over time (until convergence has been achieved), as the difference between the simulated model outputs, $f(\bm{\uptheta})$, and measurement data, $\widetilde{\textbf{d}}$, is expected to decrease with the evolution of the Markov chain.

As mentioned above, the deterministic component, $\textbf{K} \textbf{r}_{(t-1)}$, of the jump vector, $\Delta \bm{\uptheta}_{(t-1)}$, will guide the candidate states to the ``true" parameters, $\bm{\uptheta}^{*}$, of the data-generating process. This will help shorten burn-in, yet introduce asymmetry in the sampled candidate states. We must account for this asymmetry of the Kalman-inspired proposal distribution to preserve the unique stationary distribution of the Markov chain. Equation (\ref{MH}) will help to remedy the dissimilar forward $q(\bm{\uptheta}_\text{p}|\bm{\uptheta}_{(t-1)})$ and backward $q(\bm{\uptheta}_{(t-1)}|\bm{\uptheta}_\text{p})$ jump probabilities of equation (\ref{Kalman_jump}). Alternatively, we can follow the suggestion made by \citeA{vrugt2013} to randomize the direction of the jump vector, $\Delta \bm{\uptheta}_{(t-1)}$, as follows
\begin{linenomath*}
\begin{equation}
\begin{split}
\bm{\uptheta}_\text{p} & = \bm{\uptheta}_{(t-1)} \pm
(\textbf{K} \textbf{r}_{(t-1)} + \textbf{K} \bm{\upepsilon}_{(t-1)})\\
& = \bm{\uptheta}_{(t-1)} \pm \Delta \bm{\uptheta}_{(t-1)}.
\end{split}
\label{Kalman_jump_rev}
\end{equation}
\end{linenomath*}
This modification enforces symmetry of the Kalman-inspired proposal distribution with equal selection probability of the direction of the jump vector, $\Delta \bm{\uptheta}_{(t-1)}$.

The two approaches described above are easy to implement in practice and will guarantee detailed balance of the sampled Markov chain. These so-called reversibility patches do have an undesired side-effect, i.e., they deteriorate considerably the sampling efficiency of the Kalman-inspired proposal distribution, although they are theoretically appealing. This is easily demonstrated with numerical experiments. For example, if we use equation (\ref{MH}) to account for the asymmetry of the Kalman-inspired proposal distribution, then trial moves from $\bm{\uptheta}_{(t-1)}$ to $\bm{\uptheta}_\text{p}$ will be punished heavily as 
$q(\bm{\uptheta}_\text{p}|\bm{\uptheta}_{(t-1)}) \gg q(\bm{\uptheta}_{(t-1)}|\bm{\uptheta}_\text{p})$. The Hastings correction will therefore decrease the acceptance rate of candidate states. Unfortunately, the symmetric Kalman kernel of equation (\ref{Kalman_jump_rev}) provides no solace. Many of the trial moves of this proposal distribution will go to waste (especially when $\bm{\uptheta}_\text{p} = \bm{\uptheta}_{(t-1)} - \Delta \bm{\uptheta}_{(t-1)}$) as they lead the chain away from the target distribution. Thus, both reversibility patches may defeat the purpose of the Kalman-inspired proposal distribution. To take optimal advantage of the Kalman-inspired proposal distribution, we limit its application to only the first $T_\text{K}$ steps of the Markov chain. During this prescribed, relatively short burn-in period, the Kalman-inspired proposal will guide the chain to the target distribution and the resulting samples are discarded. At the end of this period we will switch the chain to a reversible proposal distribution that uses only information from past and/or present chain states to generate trial moves.

The Kalman-inspired proposal distribution can be conveniently embedded in any adequate MCMC method. In the next section, we present the implementation of this proposal distribution in the DREAM$_\text{(ZS)}$ algorithm, and an algorithmic outline of the resulting DREAM$_\text{(KZS)}$ algorithm.

\subsection{The DREAM$_\text{(KZS)}$ Algorithm}
The DREAM$_\text{(KZS)}$ algorithm is an adaptive, multiple chain MCMC method, which exploits information from past sampled states of model parameters and outputs to rapidly explore the parameter space in pursuit of the target distribution, $p(\bm{\uptheta}|\widetilde{\textbf{d}})$. This method is an extension of the DREAM$_\text{(ZS)}$ algorithm, which in turn has its roots within DE-MC$_\text{(ZS)}$ of \citeA{terbraak2008}. The next two subsections introduce the DREAM$_\text{(KZS)}$ algorithm and provide an algorithmic recipe of this new MCMC sampler. 

\subsubsection{DREAM$_\text{(KZS)}$: In Words}
Let $\bm{\Theta}_{(1)} = \{\bm{\uptheta}^{1}_{(1)},\ldots,\bm{\uptheta}^{N}_{(1)} \}$ be a $N \times k$ matrix with the initial states of the $N$ chains. These states are drawn from the prior distribution. Similarly, let $\textbf{Z} = \{ \bm{\uptheta}^{1} , \ldots, \bm{\uptheta}^{m_{0}}\}$ be an archive with $m_{0}$ ($m_{0} \gg N$) draws from the prior distribution. If $\mathcal{A}$ is a subset of $k^{-}$ dimensions of the original parameter space, $\mathbb{R}^{k^{-}}\subseteq \mathbb{R}^k$, then a candidate, $\bm{\uptheta}_\text{p}^{i}$, in the $i$th chain, $i = 1,\ldots,N$, at iteration $t$, $t = 2,\ldots,T$, is calculated from the samples in the archive, $\textbf{Z}$, using a mix of parallel direction \cite{price2006,storn1997,vrugt2008b,vrugt2009}
\begin{linenomath*}
	\begin{equation}
	\begin{array}{l}
	\Delta \bm{\uptheta}^{i}_{l \in\mathcal{A}} = \bm{\zeta}_{l} + (1 + \bm{\lambda}_{l}) \gamma_{\text{P}(\delta,k^{-})} \sum_{j=1}^{\delta} \bigl( \textbf{Z}_{l}^{\textbf{a}_{j}} - \textbf{Z}_{l}^{\textbf{b}_{j}} \bigr) \\
	\Delta \bm{\uptheta}^{i}_{l \notin\mathcal{A}} = 0 
	\end{array}
	\quad \mbox{if } u \in [0,p_\text{P}],
	\label{PDJ_DREAM_KZS}
	\end{equation}
\end{linenomath*}
snooker \cite{laloy2012,terbraak2008,vrugt2016}
\begin{linenomath*}
\begin{equation}
\Delta \bm{\uptheta}^{i} = \bm{\zeta} + (1 + \bm{\lambda})\gamma_\text{S}( \textbf{Z}^{a_{1}}_{\perp} - \textbf{Z}^{b_{1}}_{\perp}) \quad \mbox{if } u \in (p_\text{P},p_\text{P}+p_\text{S}],
\label{SJ_DREAM_KZS}
\end{equation}
\end{linenomath*}
and Kalman trial moves
\begin{linenomath*}
\begin{equation}
\Delta \bm{\uptheta}^{i} = \gamma_\text{K}\big[\textbf{K}_{(t-1)} \textbf{r}_{(t-1)}^{i} + \textbf{K}_{(t-1)} \bm{\upepsilon}_{(t-1)}^{i}\big] \quad \mbox{if } u \in (p_\text{P}+p_\text{S},1] \text{ and } t \leq T_\text{K},
\label{KJ_DREAM_KZS}
\end{equation}
\end{linenomath*}
where $u$ is a random draw from a standard uniform distribution, 
$u \sim \mathcal{U}(0,1)$; $p_\text{P}$ and $\gamma_\text{P} = 2.38/\sqrt{2\delta k^{-}}$, $p_\text{S}$ and $\gamma_\text{S} \sim \mathcal{U}(1.2,2.2)$, and $p_\text{K} = 1 - p_\text{P} - p_\text{S}$ and $\gamma_\text{K} = 1$, signify the selection probabilities and jump rates of the parallel direction, snooker and Kalman-inspired proposal distributions, respectively; $\delta$ denotes the number of chain pairs of the parallel direction proposal; $\textbf{a} = \{ a_{1},\ldots,a_{\delta} \}$ and $\textbf{b} =\{ b_{1},\ldots,b_{\delta} \}$ are $\delta$-vectors with integers randomly drawn from $\{1,\ldots,m\}$ without replacement; $\bm{\lambda}$ and $\bm{\zeta}$ are sampled independently from $\mathcal{U}_{k}(-\nu,\nu)$ and $\mathcal{N}_{k}(0,\beta)$, respectively, with $\nu = 0.05$ and $\beta = 10^{-12}$ are small compared to the width of the target distribution; $\textbf{Z}^{a_{1}}_{\perp}$ and $\textbf{Z}^{b_{1}}_{\perp}$ are orthogonal projection points of the samples $\textbf{Z}^{a_{1}}$ and $\textbf{Z}^{b_{1}}$ onto the line going through the current state of the $i$th chain, $\bm{\uptheta}^{i}_{(t-1)}$, and sample $\textbf{Z}^{c}$ of the external archive, where $c$ is an integer randomly drawn from $\{1,\ldots,m\}$ and $c \neq a_{1} \neq b_{1}$; $\textbf{r}_{(t-1)}^{i}$ and $\bm{\upepsilon}_{(t-1)}^{i}$ are residual vector of parameters $\bm{\uptheta}^{i}_{(t-1)}$ and a random draw from the distribution of measurement errors. Furthermore, to enhance the probability of a direct jump between disconnected modes of the target distribution, we use $\gamma_\text{P} = 1$ in 20\% of the parallel direction proposals. 

The candidate state of chain $i$ at iteration $t$ then becomes
\begin{linenomath*}
\begin{equation}
\bm{\uptheta}^{i}_\text{p} = \bm{\uptheta}^{i}_{(t-1)} + \Delta \bm{\uptheta}^{i},
\label{DREAM_KZS-candidate}
\end{equation}
\end{linenomath*}
and the Metropolis ratio of equation (\ref{RWM}) is used to determine whether to accept this candidate or not. If the acceptance probability, $p_{\rm acc}(\bm{\uptheta}^{i}_{(t-1)} \to \bm{\uptheta}^{i}_\text{p})$, is larger than or equal to a uniform random label, $u \sim \mathcal{U}(0,1)$, then the candidate state is accepted and the $i$th chain moves to the new position, i.e., $\bm{\uptheta}^{i}_{(t)} = \bm{\uptheta}^{i}_\text{p}$, otherwise, $\bm{\uptheta}^{i}_{(t)} = \bm{\uptheta}^{i}_{(t-1)}$. The snooker candidates demand a multiplicative, Hastings-type correction, $\alpha(\bm{\uptheta}^{i}_{(t-1)} \to \bm{\uptheta}^{i}_\text{p})$, of the acceptance probability
\begin{linenomath*}
\begin{equation}
\alpha(\bm{\uptheta}^{i}_{(t-1)} \to \bm{\uptheta}^{i}_\text{p}) = \left( \frac{ \lVert \bm{\uptheta}^{i}_\text{p} - \textbf{Z}^{c} \rVert } { \lVert \bm{\uptheta}^{i}_{(t-1)} - \textbf{Z}^{c} \rVert } \right)^{(k - 1)},
\label{snooker-correction}
\end{equation}
\end{linenomath*}
where $\lVert \cdot \rVert$ signifies the Euclidean operator. The use of the corrected acceptance probability, $\alpha(\bm{\uptheta}^{i}_{(t-1)} \to \bm{\uptheta}^{i}_\text{p})p_{\rm acc}(\bm{\uptheta}^{i}_{(t-1)} \to \bm{\uptheta}^{i}_\text{p})$, will negate bias in the jump direction of the snooker move and guarantee reversibility of the sampled chains \cite{terbraak2008}. 

The convergence speed of DREAM$_\text{(KZS)}$ is largely determined by the samples in the archive $\textbf{Z}$. 
Each $K$ generations ($K \gg 1$), we augment this archive of past states with the current states, $\bm{\Theta}_{(t)} = \{ \bm{\uptheta}^{1}_{(t)},\ldots,\bm{\uptheta}^{N}_{(t)} \}$, of the $N$ chains. This turns the DREAM$_\text{(KZS)}$ algorithm into an adaptive MCMC method with scale and orientation of the proposal that depend on the cumulative search history of the sampled chains. This proposal adaptation violates the Markovian properties, i.e., the next state of the chain should depend only on the current one. Nevertheless, adaptation has been employed in various MCMC algorithms \cite[just name a few]{haario2001,haario2006,vrugt2003}, to enhance search efficiency and protect against an inadequate selection of the proposal distribution. These algorithms are in theory not Markovian either. In practice, when the target distribution is bounded and the adaptation is performed under some regularity conditions, the adaptive MCMC algorithms can still retain the desired stationary distribution. In DREAM$_\text{(KZS)}$, if $m$ denotes the current number of samples in the archive, then appending $N$ states changes $\textbf{Z}$ by an order of $K/t$, which decreases with iteration $t$. The three proposal distributions used in DREAM$_\text{(KZS)}$ thus become invariant as the length of the thinned past increases without bound. This so-called diminishing adaptation ensures ergodic chains that converge to the exact target distribution \cite{roberts2007}. 

The $k^{-}$ members of the subset $\mathcal{A}$ are drawn at random from the entries $\{1,\ldots,k\}$ with the help of a crossover probability, $\eta \in (0,1]$. Each time a proposal is generated, we draw a $k$-vector, $\textbf{u} = \{u_{1},\ldots,u_{k}\}$, of standard uniform labels, $\textbf{u} \sim \mathcal{U}_{k}(0,1)$. All entries $j$ of $\textbf{u}$ which satisfy $u_{j} \leq \eta$ are stored in the subset $\mathcal{A}$ and used to span the search subspace that will be sampled using equation (\ref{PDJ_DREAM_KZS}). $\mathcal{A}$ must at least have one element, otherwise the jump vector will have zero length. To enhance search efficiency, we use a geometric series of $n_{\eta}$ crossover values, $\bm{\upeta} = \{ 1/n_{\eta},2/n_{\eta},\ldots,1\}$, and sample the value of $\eta$ from a discrete multinomial distribution, $\mathcal{M}(\bm{\upeta},\textbf{p}_{\eta})$, on $\bm{\upeta}$ with selection probabilities $\textbf{p}_{\eta}$. The values of $\textbf{p}_{\eta}$ are tuned adaptively during burn-in by maximizing the traveled distance of the $N$ chains. This adaptation scheme is described in detail in \citeA{vrugt2008a} and \citeA{vrugt2009}. The use of a vector of crossover probabilities enables single-site Metropolis ($\mathcal{A}$ has one element), Metropolis-within-Gibbs ($\mathcal{A}$ has one or more elements) and regular Metropolis sampling ($\mathcal{A}$ has $k$ elements), and enables chains to sample outside the subspace spanned by their current positions. The default setting of $n_{\eta} = 3$ has shown to work well in practice.

We are now left with a numerical implementation of the Kalman-inspired proposal distribution. We use this proposal distribution only during the first $T_\text{K}$ generations, after which $p_\text{K} = 0$ and the selection probabilities of the parallel direction and snooker proposal distributions, $p_\text{P}$ and $p_\text{S}$, respectively, are re-normalized so that their values add up to one. The Kalman gain, $\textbf{K}_{(t-1)}$, in equation (\ref{KJ_DREAM_KZS}) is computed from an archive $\textbf{Z}_\text{K}$ with $m_\text{K}$ samples drawn from the chain history, as follows 
\begin{linenomath*}
\begin{equation}
\textbf{K}_{(t-1)}=\text{Cov}\bigr(\textbf{Z}_\text{K},f(\textbf{Z}_\text{K})\bigr)\bigr[\text{Cov}\bigr(f(\textbf{Z}_\text{K}),f(\textbf{Z}_\text{K})\bigr)+\textbf{R}\bigr]^{-1},
\label{KJ_implementation}
\end{equation}
\end{linenomath*}
where $f(\textbf{Z}_\text{K})$ is a $m_\text{K}\times n$ matrix with model outputs of the $m_\text{K}$ samples in the archive $\textbf{Z}_\text{K}$.

It should be noted here that although the Kalman filter is based on Gaussian assumption, an MCMC algorithm (e.g., DREAM in this work) that adopts the Kalman-inspired proposal distribution can be used to explore complex and non-Gaussian (e.g., multi-modal) target distributions (see section 3.4). Due to the fact that they are both based on the Kalman formula, one might confuse the Kalman-inspired proposal distribution that only generates a candidate at a time, with  ensemble Kalman filter (EnKF) and its variants \cite{evensen2009data} that use an ensemble of parameter samples updated with the Kalman formula to directly approximate the posterior. In EnKF and its variants, one should use the prior covariances to calculate the Kalman gain, then update the prior ensemble with the Kalman formula, and finally use the updated ensemble to approximate the mean and covariance of the posterior. These methods are restricted to problems with Gaussian parameter distributions. However, in our work, the Kalman-inspired proposal is not designed to directly approximate the posterior, but to suggest where the Markov chain might move at the next step by generating a candidate. At the very beginning of the MCMC simulation, it is very likely that the candidate states are far from the posterior mean. Implemented in an adaptive Metropolis sampler, the Kalman-inspired proposal actually uses covariances that undergo adaptation. The covariances can converge, after burn-in, towards posterior covariances. To our best knowledge, proposal distributions used in various popular MCMC algorithms also adopt Gaussian forms, which does not mean that the MCMC algorithms are restricted to problems with normal posterior distributions. 
	
Nevertheless, in complicated, multi-modal problems, the Kalman-inspired proposal faces the risk of missing secondary modes. To improve the capability of the Kalman-inspired proposal in multi-modal cases, one can modify the original Kalman-inspired proposal by adopting the local Kalman update strategy proposed by \citeA{zhang2018iterative}, which is specifically designed to solve inverse problems with (possible) multi-modal parameter distributions. The basic idea behind the local Kalman update strategy is simple: Although globally the parameter distribution might be non-Gaussian, or even multi-modal, one can still use a Gaussian distribution to describe the local parameter distribution. For example, when overall the parameter distribution is multi-modal, if one only looks at the neighborhood of a single mode, the local parameter distribution can be still close to Gaussian. Using an integrated measure of distance proposed by \citeA{zhang2018iterative}, one can find the local ensemble of each sample in the archive for the Kalman-inspired proposal. Through updating all the local ensembles separately, and generating a candidate randomly from the multiple updated local ensembles, the modified proposal distribution will have a larger chance to sample the secondary modes. As the original Kalman-inspired proposal distribution already works well in multi-modal cases (see section 3.4), to prevent the length of the paper from being too long, the local Kalman update strategy is not tested here. We refer the reader to \citeA{zhang2018iterative} for a detailed description of the related method. Nevertheless, in the Kalman-inspired proposal, distribution of measurement errors still needs to be Gaussian. To enable a proper use of the Kalman-inspired proposal, one can transform (possible) non-Gaussian measurement errors to Gaussian variables with some transformation method \cite{chou1998,sakia1992}.

\subsubsection{DREAM$_\text{(KZS)}$: Algorithmic Recipe}
This section provides an algorithmic recipe of the DREAM$_\text{(KZS)}$ algorithm (see Algorithm \ref{alg:DREAM_KZS}). This recipe translates the initial population, $\bm{\Theta}_{(1)} = \{ \bm{\uptheta}^{1}_{(1)},\ldots,\bm{\uptheta}^{N}_{(1)}\}$, into samples from the target distribution. 

The parallel direction, snooker and Kalman-inspired proposal distributions are used interchangeably during the first $T_\text{K}$ generations to propose candidate states in the $N$ chains. After this prescribed burn-in period, the selection probability of the Kalman-inspired proposal distribution is set to zero, i.e., $p_\text{K} = 0$, and the selection probabilities of the parallel direction and snooker proposal distributions, $p_\text{P}$ and $p_\text{S}$, are normalized to sum to unity. Thus, after the relatively short burn-in period, the sampled chains will maintain detailed balance.

Convergence of the sampled chains of the DREAM$_\text{(KZS)}$ algorithm can be monitored with a variety of different metrics, including within-chain and between-chain statistics. Of these, the univariate, $\widehat{R}_{j}, j = 1,\ldots,k$, and multivariate, $\widehat{R}^{k}$, scale reduction factors of \citeA{gelman1992} and \citeA{brooks1998}, respectively, are most widely used to assess convergence for multiple chain methods. These two diagnostics compare for each parameter individually, or for the distribution as a whole, the within-chain and between-chain variances or covariance matrices, respectively. Values of $\widehat{R}_{j} \leq 1.2, \forall j \in \{1,\ldots,k\}$ and $\widehat{R}^{k} \leq 1.2$ demonstrate convergence of the sampled chains to a stationary distribution. 

\begin{algorithm}[H]
	\caption{DREAM$_\text{(KZS)}$}
	\label{DREAM_KZS}
	\begin{algorithmic}[1]
		\STATE Define $\bm{\upeta} = \{1/n_\eta, 2/n_{\eta},\ldots,1\}$ and $\textbf{p}_{\eta} = \textbf{1}_{n_\eta}/n_{\eta}$. Set $n_\eta = 3$, $m_{0} = 10k$ and $K = 10$.
		\STATE Draw an initial archive, $\textbf{Z}$, from the prior distribution, $\bm{\uptheta}^{j} \sim p(\bm{\uptheta})$, $j = 1,\ldots,m_{0}$.
		\STATE Draw the initial states of the $N$ chains from the prior distribution, $\bm{\uptheta}^{i}_{(1)} \sim p(\bm{\uptheta})$, $i = 1,\ldots,N$.
		\STATE Initialize at zero the traveled distance, $D_{l} = 0$, of each crossover value, $l = 1,\ldots,n_\eta$.
		\FOR{$t = 2,\ldots,T$}
		\STATE Compute $\bm\Delta_{(t-1)}$, the $k$-vector of standard deviations of $\bm{\Theta}_{(t-1)} = \{\bm{\uptheta}^{1}_{(t-1)},\ldots,\bm{\uptheta}^{N}_{(t-1)}\}$.
		\STATE Draw $N$ labels, $\textbf{u} = \{u^{1},\ldots,u^{N}\}$, from a standard uniform distribution, $u^{i} \sim \mathcal{U}(0,1)$.
		\STATE Draw an integer, $J$, from the multinomial distribution, $\mathcal{M}(\{1,2,3\},\{p_\text{P},p_\text{S},p_\text{K}\})$.
		\STATE If $J$ = 1, draw a crossover value, $\eta_{l} \sim \mathcal{M}(\bm{\upeta},\textbf{p}_{\eta})$, otherwise, $l = n_{\eta}$, thus $\eta_{l} = 1$.
		\FOR{$i = 1,\ldots,N$}
		\STATE Create a candidate, $\bm{\uptheta}^{i}_{\text{p}}$, using the parallel direction ($J=1$), snooker ($J=2$) or Kalman-inspired  ($J=3$) proposal distribution.
		\STATE Calculate the Metropolis ratio, $p_{\text{acc}}(\bm{\uptheta}^{i}_{(t-1)} \to \bm{\uptheta}^{i}_\text{p})$, using equation (\ref{RWM}).
		\STATE If $J$ = 2, compute $\alpha(\bm{\uptheta}^{i}_{(t-1)} \to \bm{\uptheta}^{i}_\text{p})$ using equation (\ref{snooker-correction}), otherwise, $\alpha(\bm{\uptheta}^{i}_{(t-1)} \to \bm{\uptheta}^{i}_\text{p}) = 1$.
		\STATE If $\alpha(\bm{\uptheta}^{i}_{(t-1)} \to \bm{\uptheta}^{i}_\text{p})p_{\text{acc}}(\bm{\uptheta}^{i}_{(t-1)} \to \bm{\uptheta}^{i}_\text{p})\ge u^{i}$, set $\bm{\uptheta}^{i}_{(t)} = \bm{\uptheta}^{i}_\text{p}$, otherwise, set $\bm{\uptheta}^{i}_{(t)} = \bm{\uptheta}^{i}_{(t-1)}$.
		\ENDFOR
		\STATE Update the traveled distance, $D_{l} = D_{l} + \sum_{i=1}^{N} || (\bm{\uptheta}^{i}_{(t)} - \bm{\uptheta}^{i}_{(t-1)})/\bm{\Delta}_{(t-1)} ||$ for $\eta_{l}$. 
		\STATE If $\text{mod}(t,K) = 0$, append the current chain states, $\bm{\Theta}_{(t)} = \{ \bm{\uptheta}^{1}_{(t)},\ldots,\bm{\uptheta}^{N}_{(t)}\}$, to $\textbf{Z}$ and set $m = m + N$. 
		\STATE If $t \leq T_\text{K}$, update $\textbf{p}_{\eta}$ using the cumulative traveled distance, $\textbf{D} = \{D_{1},\ldots,D_{n_\eta}\}$, of $\bm{\upeta}$, otherwise, set $p_\text{K} = 0$, $p_\text{T} = p_\text{P} + p_\text{S}$, $p_\text{P} = p_\text{P}/p_\text{T}$ and $p_\text{S} = p_\text{S}/p_\text{T}$.
		\ENDFOR
	\end{algorithmic}
    \label{alg:DREAM_KZS}
\end{algorithm}{\tiny }

\section{Illustrative Case Studies}
In this section, we analyze, test and evaluate the performance of the DREAM$_\text{(KZS)}$ algorithm by application to several different case studies with complex, high-dimensional or multi-modal target distributions. These studies involve mathematical models of rainfall-runoff transformation, groundwater flow and contaminant transport. We use simulated data of known parameter values to evaluate the performance of the DREAM$_\text{(KZS)}$ algorithm and benchmark its sampling efficiency against the DREAM$_\text{(ZS)}$ algorithm. The DREAM$_\text{(ZS)}$ algorithm has found widespread application and use, and many published studies corroborate its excellent ability to rapidly sample complex, high-dimensional target distributions \cite{bikowski2012,muleta2012,ramin2014,shi2014,wohling2011,zhang2018}. In the numerical experiments, we suggest to set $p_\text{P} = 0.6$, $p_\text{S} = 0.1$, $p_\text{K} = 0.3$, and $T_\text{K} = 0.3T$ in the DREAM$_\text{(KZS)}$ algorithm, where $T$ signifies the maximum number of chain generations. Moreover, to demonstrate that the Kalman-inspired proposal distribution can be conveniently embedded in any adequate MCMC method, in section 3.1, we further test the performance of the Kalman-inspired proposal in two plainer MCMC algorithms, i.e., AM \cite{haario2001} and Delayed Rejection AM \cite<DRAM,>[]{haario2006}.

\subsection{Case Study 1: A Rainfall-Runoff Model}
Our first case study considers application of the DREAM$_\text{(KZS)}$ algorithm to modeling of the rainfall-discharge relationship of the Guadalupe River basin \cite{duan2006}. We use the seven-parameter hmodel of \citeA{schoups2010} to simulate daily discharge records of the Guadalupe River at Spring Branch, Texas using basin average estimates of precipitation and potential evapotranspiration. The hmodel transforms rainfall into runoff at the watershed outlet using four different control volumes and simulation of interception, throughfall, evaporation, runoff generation, percolation through surface and subsurface routings. The seven parameters of the hmodel and their prior ranges are listed in Table \ref{table:1}. We refer the reader to \citeA{schoups2010} for a detailed description of the hmodel, including model structure, process specification and parameterization.
\vspace*{0px}
\begin{center} [ Table 1 NEAR HERE ] \end{center}
\vspace{0px}

We use the hmodel parameter values listed in Table \ref{table:1} to simulate a $n = 1,827$ record of daily discharge values, $\textbf{d}$. This record is corrupted with heteroscedastic measurement errors by drawing from the $n$-variate normal distribution, $\widetilde{\textbf{d}} \sim \mathcal{N}_{n}(\textbf{d},\textbf{R})$, with values of $(0.05\textbf{d})^{2}$ on the main diagonal of the covariance matrix, $\textbf{R}$, and zero entries elsewhere. We now use the perturbed record, $\widetilde{\textbf{d}}$, to infer the posterior distribution of the hmodel parameters with the DREAM$_\text{(KZS)}$ algorithm using $N = 4$ chains with $T = 6,000$ samples in each chain. We assume a uniform prior parameter distribution over the ranges listed in Table \ref{table:1}, and use the Gaussian likelihood function of equation (\ref{likelihood}) with known measurement error covariance matrix, $\textbf{R}$.

\begin{figure}
	\centering
	\includegraphics[width=35pc]{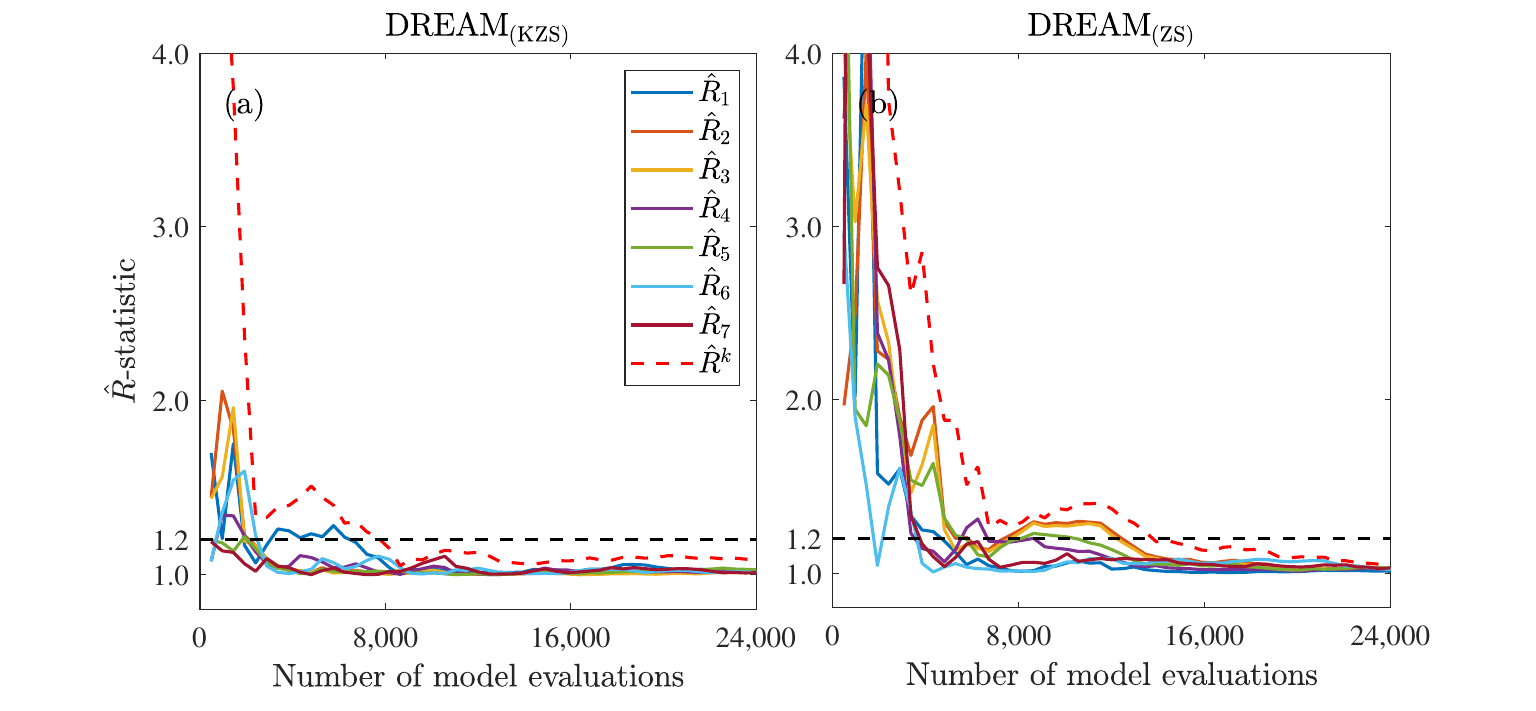}
	\caption{Traces of the $\widehat{R}$-statistics of the seven hmodel parameters using the (a) DREAM$_\text{(KZS)}$ and (b) DREAM$_\text{(ZS)}$ algorithms. The horizontal black dashed lines depict the threshold value of 1.2 used to declare convergence to a stationary distribution.}
	\label{fig:1}
\end{figure}

Figure \ref{fig:1} presents a trace plot of the univariate convergence diagnostic for each of the $k = 7$ hmodel parameters, $\hat{R}_{j}, j=1,...,k$, and the multivariate convergence diagnostic, $\widehat{R}^{k}$, using the DREAM$_\text{(KZS)}$ (left panel) and DREAM$_\text{(ZS)}$ (right panel) algorithms. The different hmodel parameters are coded with different colors. The horizontal black dashed line in each panel demarcates the threshold of $1.2$ below which the chains are assumed to have converged to a stationary distribution. The DREAM$_\text{(KZS)}$ algorithm requires about 1,800 iterations (i.e., 7,200 model evaluations) to satisfy the stipulated convergence threshold of 1.2. The DREAM$_\text{(ZS)}$ algorithm, on the contrary, needs a substantially larger number of about 3,360 iterations (i.e., 13,440 model evaluations) to converge to the stationary distribution. These findings suggest an improvement in simulation efficiency on the order of 1.87 times. 
\vspace*{0px}
\begin{center} [ Figure 1 NEAR HERE ] \end{center}
\vspace{0px}

\begin{figure}
	\centering
	\includegraphics[width=35pc]{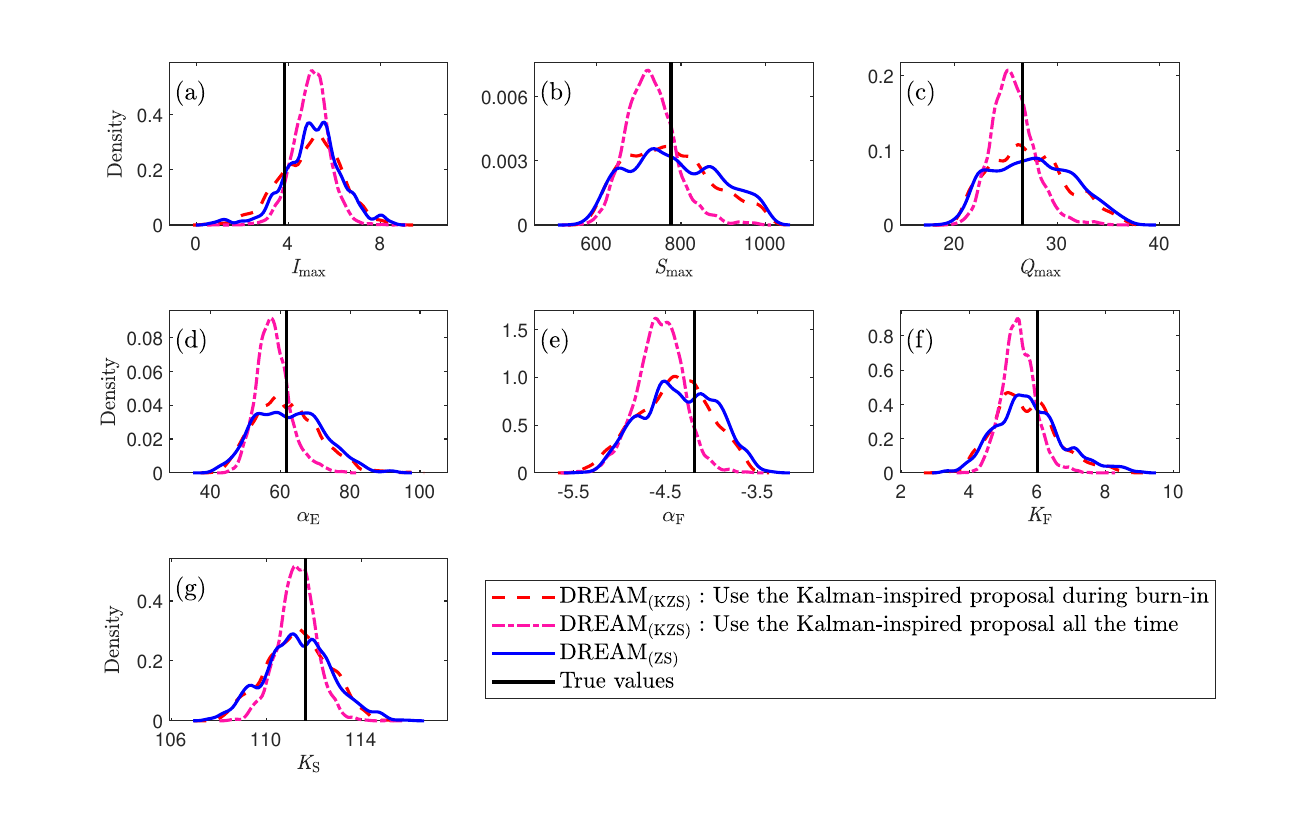}
	\caption{Marginal posterior distributions of the seven hmodel parameters obtained by DREAM$_\text{(KZS)}$ that uses the Kalman-inspired proposal during burn-in (default setting, red dashed lines), DREAM$_\text{(KZS)}$ that uses the Kalman-inspired proposal all the time (magenta dash-dotted lines), and DREAM$_\text{(ZS)}$ (reference results, blue lines), respectively. The vertical black lines display the ``true" values of the hmodel parameters.}
	\label{fig:2}
\end{figure}

These initial results are encouraging, yet they would not mean much if the DREAM$_\text{(KZS)}$ algorithm did not approximate correctly the target distribution. In Figure \ref{fig:2}, we plot marginal posterior distributions of the seven hmodel parameters derived from the DREAM$_\text{(KZS)}$ (red dashed lines) and DREAM$_\text{(ZS)}$ (blue lines) algorithms. The densities are estimated with a normal kernel function using the last $1,000$ samples in each of the $N = 4$ Markov chains. The ``true" values of the hmodel parameters are separately indicated in each panel with a vertical black line. The approximated marginal distributions of both algorithms appear remarkably similar and center nicely on the ``true" values of the hmodel parameters. This provides evidence for the claim that the DREAM$_\text{(KZS)}$ algorithm successfully approximates the target distribution.
\vspace*{0px}
\begin{center} [ Figure \ref{fig:2} NEAR HERE ] \end{center}
\vspace{0px}

Now one question still remains, i.e., what if we do not switch the Kalman-inspired proposal distribution off after burn-in? To answer this question, we further implement the DREAM$_\text{(KZS)}$ algorithm that uses the Kalman-inspired proposal all the time. From Figure S1 (in the supporting information) we can find that, the $\widehat{R}$-diagnostic plots for the two DREAM$_\text{(KZS)}$ approaches can reach the threshold value of 1.2 with a similar number of iterations. This is not surprising as the two approaches use the same settings during burn-in. After burn-in, if we continue to use the Kalman-inspired proposal, we can obtain slightly smaller $\widehat{R}$ values than the approach that switches off the Kalman-inspired proposal after burn-in. However, as shown in Figure \ref{fig:2}, using the Kalman-inspired proposal all the time in MCMC will make the final estimation of posterior (magenta dash-dotted lines) slightly narrower than the reference results obtained by DREAM$_\text{(ZS)}$ (blue lines), while the recommended approach that uses the Kalman-inspired proposal only during burn-in obtains reliable results (red dashed lines). To investigate the performance of the Kalman-inspired proposal, we further analyze the acceptance rates of the two DREAM$_\text{(KZS)}$ approaches. In Figure S2, we depict the evolution of the numbers of accepted candidate states in the two DREAM$_\text{(KZS)}$ approaches. During the prescribed burn-in period (i.e., the first 30\% of the MCMC simulation), the average acceptance rate (using the parallel direction, snooker and Kalman-inspired proposals with suggested selection probabilities) is about 0.17; if we only account for the candidate states generated by the Kalman-inspired proposal, the acceptance rate can be as high as about 0.40, which indicates that the Kalman-inspired proposal can generate ``good” candidates that are less likely to be rejected. After burn-in, if we switch off the Kalman-inspired proposal, the average acceptance rate (using the other two proposals) during the last 70\% of the MCMC simulation will be about 0.10; if the Kalman-inspired proposal is still left on, the average acceptance rate (using the three proposals) can be about 0.33. The above results indicate that the Kalman-inspired proposal can shorten burn-in, but it introduces asymmetry to the sampled chains. To maintain detailed balance of the Markov chains, we can adopt the three strategies proposed in section 2.2, among which the simplest way is to restrict the Kalman-inspired proposal distribution to the burn-in period, which is our default setting in the following tests.

\begin{figure}
	\centering
	\includegraphics[width=36pc]{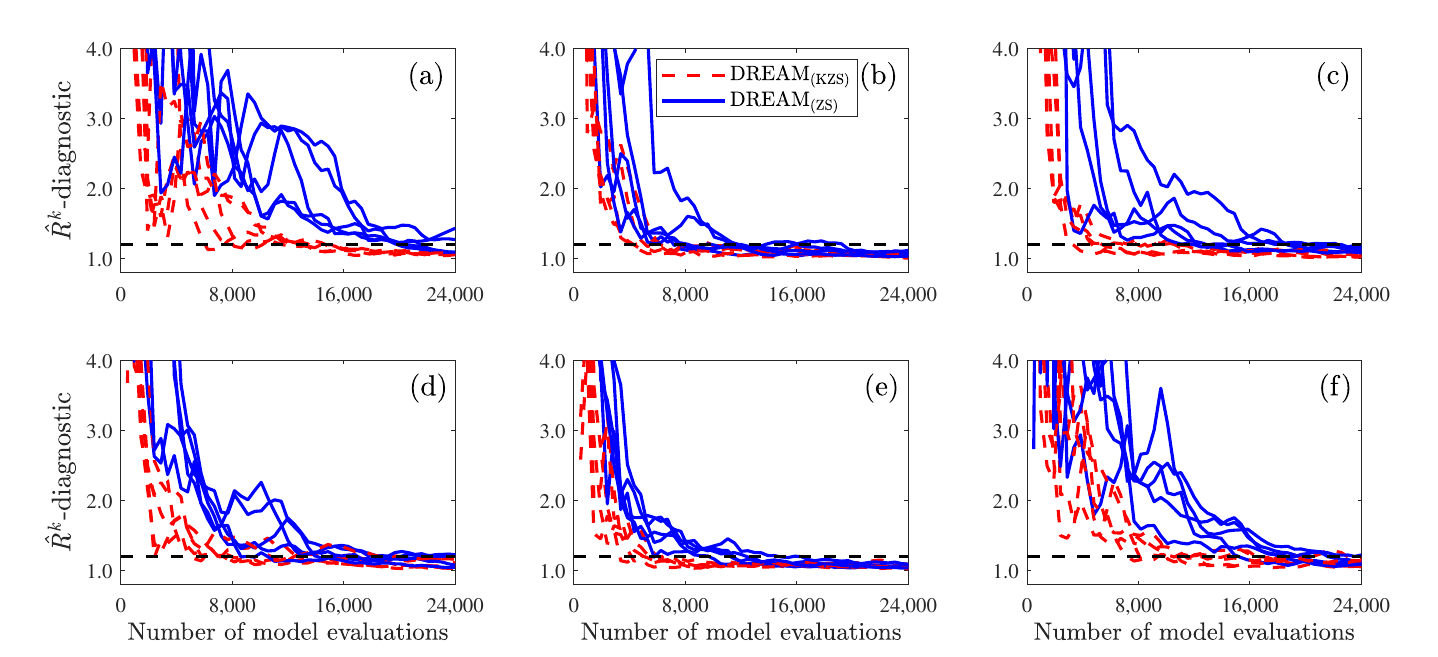}
	\caption{Evolution of the multivariate $\widehat{R}^k$-convergence diagnostic derived from the DREAM$_\text{(KZS)}$ (red dashed lines) and DREAM$_\text{(ZS)}$ (blue lines) algorithms for six discharge data sets (between plots) each with five repetitions (within plot). The horizontal black dashed line in each plot displays the threshold of 1.2 used to declare convergence to a stationary distribution.}
	\label{fig:3}
\end{figure}

To generalize our findings, we repeat our numerical experiment using six sets of other hmodel parameter vectors drawn randomly from the uniform prior distribution. To negate sampling variability, we repeat the simulation of the DREAM$_\text{(KZS)}$ and DREAM$_\text{(ZS)}$ algorithms five times for each parameter vector. The results of our analysis are presented in Figure \ref{fig:3}, which presents traces of the $\hat{R}^{k}$-diagnostic for the five repetitions (within a graph) of each parameter vector (between graphs) using the DREAM$_\text{(KZS)}$ and DREAM$_\text{(ZS)}$ algorithms. Thus, each subplot corresponds to a different hmodel parameter vector and contains a separate trace of the $\hat{R}^{k}$-convergence diagnostic for each of its five repetitions with the same set of measurement data. Color coding differentiates between the DREAM$_\text{(KZS)}$ (red dashed lines) and DREAM$_\text{(ZS)}$ (blue lines) algorithms. As we plot the multivariate $\hat{R}^{k}$-statistic, a single line suffices for each combination of algorithm, parameter vector and repetition simulation. The results in Figure \ref{fig:3} generalize our earlier findings. Indeed, the DREAM$_\text{(KZS)}$ algorithm consistently requires fewer model evaluations to converge to the posterior distribution. Note, that the speed-up in sampling efficiency of the DREAM$_\text{(KZS)}$ algorithm is hardly impressive. This is not surprising due to the rather low dimensionality of the target distribution.
\vspace*{0px}
\begin{center} [ Figure 3 NEAR HERE ] \end{center}
\vspace{0px}

\begin{figure}
	\centering
	\includegraphics[width=35pc]{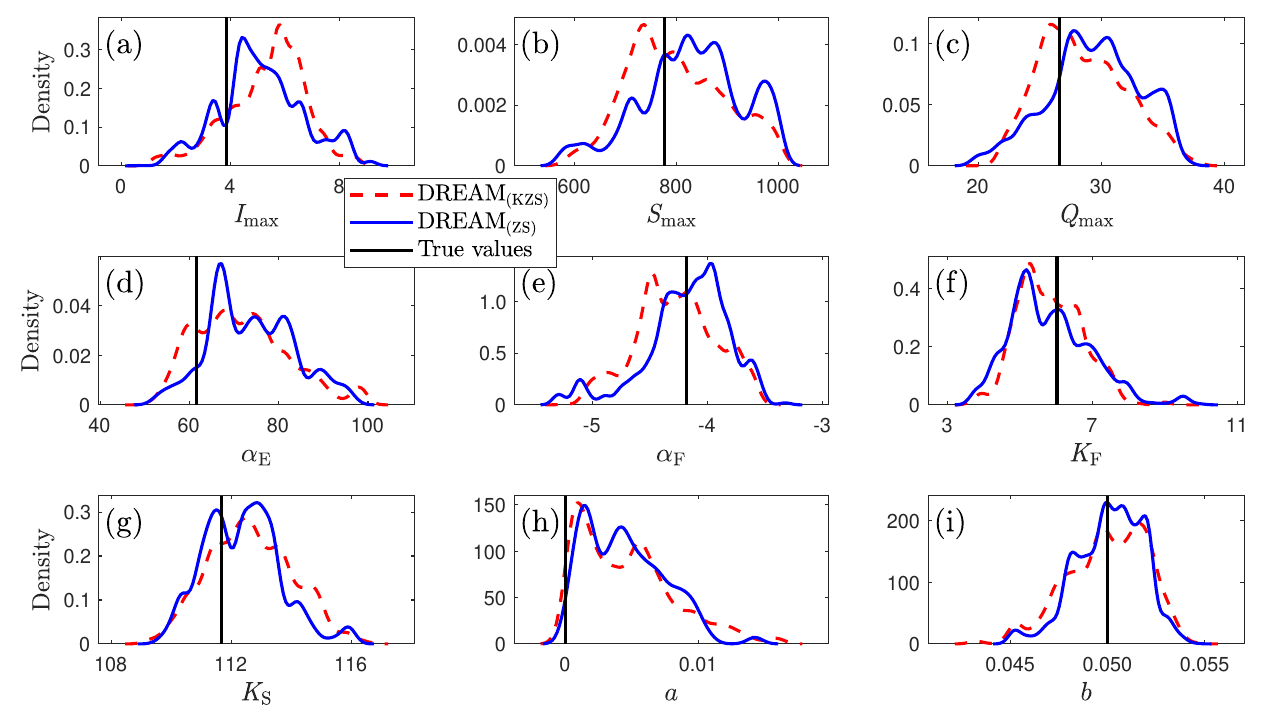}
	\caption{Marginal posterior distributions of the seven hmodel parameters and two measurement error model coefficients derived from the DREAM$_\text{(KZS)}$ (red dashed lines) and DREAM$_\text{(ZS)}$ (blue lines) algorithms. The ``true" values of the $k = 9$ parameters are separately indicated in each panel with a vertical black line.}
	\label{fig:4}
\end{figure}

Up to now, we have conveniently assumed the covariance matrix of measurement errors to be known, i.e., $\bm{\upepsilon} \sim \mathcal{N}_{n}(\textbf{0},\textbf{R})$, where $\sqrt{\textbf{R}}$ was diagonal with entries equal to 1/20 of the ``true" model responses, $\textbf{d}$. In practice, the distribution of measurement errors may not be known \textit{a-priori}. In those cases, we can define a measurement error model, say, the standard deviation of the measurement errors is expressed as, $\bm{\upsigma} = a + b\widetilde{\textbf{d}}$. The coefficients $a$ and $b$ can be treated as nuisance variables whose values can be inferred simultaneously with the model parameters \cite<see e.g., >[]{schoups2010,vrugt2016}. As a proof of concept, Figure \ref{fig:4} presents marginal posterior distributions of the hmodel parameters and coefficients $a$ and $b$ derived from the DREAM$_\text{(KZS)}$ (red dashed lines) and DREAM$_\text{(ZS)}$ (blue lines) algorithms. The estimated marginal distributions of both algorithms are in close agreement with each other and cover the assumed values (vertical black lines) of the hmodel parameters and measurement error model coefficients ($a = 0$ and $b = 0.05$).
\vspace*{0px}
\begin{center} [ Figure \ref{fig:4} NEAR HERE ] \end{center}
\vspace{0px}

It is noted here that the Kalman-inspired proposal distribution is not tied to a particular MCMC algorithm, e.g., DREAM, but is straightforwardly extendable to any adequate MCMC method. To support this claim, we also introduce the Kalman-inspired proposal distribution to two widely used, single chain MCMC methods, i.e., the AM and DRAM algorithms \cite{haario2001,haario2006}. When implementing the AM algorithm, we first start the Markov chain from the prior mean values and set the chain length as $T=24,000$. Introducing the Kalman-inspired proposal distribution to the AM algorithm produces a modified method, which is termed the AM$_\text{(K)}$ algorithm here. In the AM$_\text{(K)}$ algorithm, the Kalman-inspired proposal distribution is used during a prescribed burn-in period (the first 30\% of the MCMC simulation) with a selection probability of $p_\text{K}=0.3$. With the same settings, we then run the DRAM and corresponding DRAM$_\text{(K)}$ algorithms to infer the posterior distribution of the seven hmodel parameters. As shown in Figure S3, conditioned on the same set of measurement data (generated from one set of random parameter vector and corrupted with a zero-mean normally distributed error), the AM (black lines), AM$_\text{(K)}$ (black dashed lines), DRAM (blue lines), DRAM$_\text{(K)}$ (blue dashed lines) and DREAM$_\text{(KZS)}$ (red dashed lines) algorithms can obtain very consistent estimates of the posterior distribution. From the traces of the seven hmodel parameters obtained by the AM (red dots in Figure S4), AM$_\text{(K)}$ (blue dots in Figure S4), DRAM (red dots in Figure S5) and DRAM$_\text{(K)}$ (blue dots in Figure S5) algorithms, it is evident that the AM$_\text{(K)}$ and DRAM$_\text{(K)}$ algorithms require fewer iterations to reach the stationary regime than their original counterparts, and the Kalman-inspired proposal can improve the performance of AM more than DRAM. After systematic evaluations, we find that the initial state, or starting point, of the Markov chain can have a significant impact on the performance of the AM algorithm. As shown in Figure S6, starting from a random state drawn from the prior distribution, the AM algorithm (red dots) performs much worse than the AM algorithm starting from the prior mean state (red dots in Figure S4). However, starting from the same random state, the AM$_\text{(K)}$ algorithm (blue dots in Figure S6) exhibits a more desirable performance than the AM algorithm. Starting from the same random initial state, the DRAM algorithm (red dots in Figure S7) performs better than the AM algorithm, and adopting the Kalman-inspired proposal in DRAM (i.e., DRAM$_\text{(K)}$, blue dots in Figure S7) can further improve the performance. It is noted here that the use of delayed rejection in DRAM (or DRAM$_\text{(K)}$) will make the actual number of model evaluations much larger than the prescribed chain length of $T=24,000$. To improve the performance of MCMC, it is common practice to first run an optimizer to obtain a ``good" initial state. Here we adopt a simplex search method developed by \citeA{lagarias1998convergence} (the build-in function ``fminsearch" in MATLAB) to obtain a head start. The optimizer explores the parameter space from the prior mean values and calls 1,406 model evaluations in total. With the optimized initial state, we run the AM algorithm and plot the traces (red dots) of the seven hmodel parameters in Figure S8. In this figure, we again present the results of the AM$_\text{(K)}$ algorithm starting from the prior mean state (blue dots). Without extra model evaluations, the AM$_\text{(K)}$ algorithm still works slightly better than the AM algorithm with an initial optimizer. From the above simulation results we can conclude that, the starting point of the Markov chain can have a considerable effect on the convergence speed of a plain MCMC algorithm like AM. Thus, finding a good starting point, e.g., using an optimization method, can improve the simulation efficiency of the MCMC algorithm. Nevertheless, without resorting to looking for a good starting point, one can also employ an effective proposal distribution, e.g., the Kalman-inspired proposal distribution formulated in this work, to generate ``good" candidates that can accelerate the movement of the chain to the target region. To enhance search efficiency, adaptation of covariance(s) based on the chain history is always performed.

\subsection{Case Study 2: Groundwater Contaminant Source Identification}
\begin{figure}
	\centering
	\includegraphics[width=25pc]{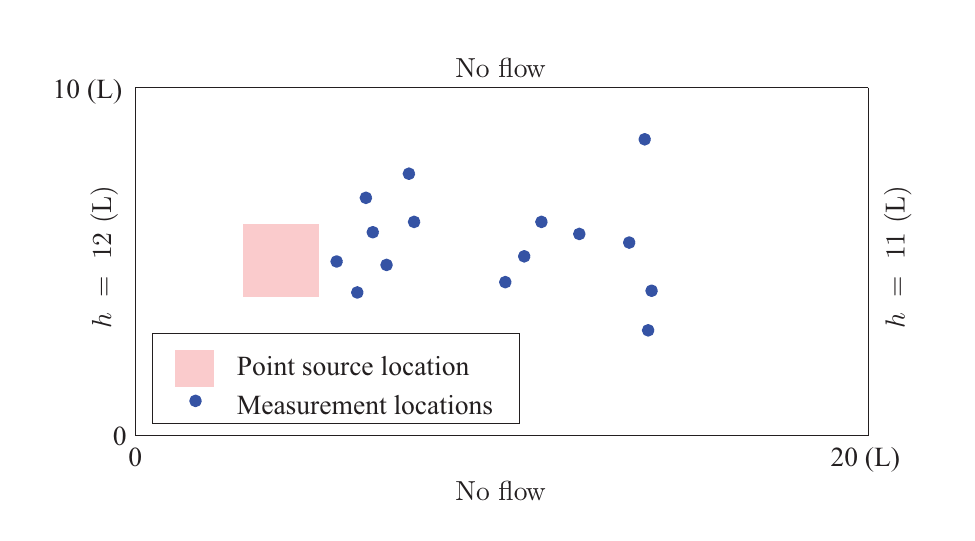}
	\caption{Schematic overview of the simulated spatial domain with imposed boundary conditions. The area with light red signifies the potential location of contaminant release from a point source, and the blue dots represent the wells at which measurements of the steady-state hydraulic head and transient contaminant concentration are collected.}
	\label{fig:5}
\end{figure}

The second case study considers two-dimensional simulation of steady-state groundwater flow and contaminant transport. We consider a hypothetical rectangular flow domain (see Figure \ref{fig:5}) with $x \in [0,20]$ (L) and $y \in [0,10]$ (L) in units of length in the horizontal and vertical direction, respectively. A constant-head of $h = 12$ (L) and $h = 11$ (L) at the left and right-side of the domain, respectively, and no-flow condition at the top and bottom of the field impose a steady-state water flux from left to right across the domain. The hydraulic conductivity field, $\textbf{K}(x,y)$ (LT$^{-1}$) of the domain is assumed to be spatially heterogeneous and isotropic with covariance function, $C_\mathcal{K}(x_1,y_1;x_2,y_2)$, of two arbitrary points, $(x_1,y_1)$ and $(x_2,y_2)$, for the log-transformed field, $\mathcal{K} = \log\bigl(\textbf{K}(x,y)\bigr)$, equal to \cite{zhang2004}
\begin{linenomath*}
	\begin{equation}
	C_\mathcal{K}(x_1,y_1;x_2,y_2) = \sigma_\mathcal{K}^2 \exp \left( -\frac{ \left| x_1 - x_2 \right| } { \lambda_x } - \frac{ \left|y_1 - y_2 \right| }{ \lambda_y} \right),
	\label{correlation}
	\end{equation}
\end{linenomath*}
where $\sigma^{2}_\mathcal{K}$ is the variance of the log-conductivity field, and $\lambda_x$ (L) and $\lambda_y$ (L) signify the correlation lengths in the $x$ and $y$ direction, respectively. 
\vspace*{0px}
\begin{center} [ Figure \ref{fig:5} NEAR HERE ] \end{center}
\vspace{0px}

We sample a reference hydraulic conductivity field (see Figure \ref{fig:8}a) with mean logarithmic conductivity, $\mu_\mathcal{K} = 2$, and values of $\sigma^{2}_\mathcal{K} = 1$, $\lambda_{x} = 10$ (L), and $\lambda_{y} = 5$ (L) in equation (\ref{correlation}). Next, we simulate the steady-state hydraulic head, $h$ (L), and corresponding pore water velocity, $v$ (LT$^{-1}$), in our rectangular domain by solving the following two equations numerically with MODFLOW \cite{harbaugh2000} 
\begin{linenomath*}
	\begin{equation}
	\frac{\partial}{\partial x_i}\left(K_i\frac{\partial h}{\partial x_i}\right) = 0,
	\label{flow1}
	\end{equation}
\end{linenomath*}
and
\begin{linenomath*}
	\begin{equation}
	v_i = -\frac{K_i}{\phi}\frac{\partial h}{\partial x_i},
	\label{flow2}
	\end{equation}
\end{linenomath*}
where $\phi = 0.25$ (-) is the aquifer porosity and the subscript $i$ signifies the respective coordinate axis ($i = 1, 2$).

To simulate contaminant transport, we pollute the water in our flow domain with an unknown point source. The exact location of this point source is assumed unknown, but its spatial coordinates, $x_\text{s}$ and $y_\text{s}$, must be found within the light red square depicted in Figure \ref{fig:5}. As the release strength of the contaminant may be time dependent, we use a simple step function to simulate transient mass-loading rates. This step function is composed of six equidistant time intervals with constant mass-loading rate, $s_{i}$ (MT$^{-1}$), $i=1,...,6$, in each segment. The space-time concentration, $C(x,y,t)$ (ML$^{-3}$), of the contaminant in our rectangular domain is simulated with MT3DMS \cite{zheng1999}, using numerical solution of the advection-dispersion equation
\begin{linenomath*}
	\begin{equation}
	\frac{\partial (\phi C)}{\partial t} = \frac{\partial}{\partial x_i} \left( \phi D_{ij}\frac{\partial C}{\partial x_j} \right) - \frac{\partial}{\partial x_i}(\phi v_iC) + q_\text{a}C_\text{s},
	\label{solute}
	\end{equation}
\end{linenomath*}
where $q_\text{a}$ (T$^{-1}$) denotes the volumetric flow rate per unit volume of the aquifer, $C_\text{s}$ (ML$^{-3}$) is the concentration of the contaminant source, and $D_{ij}$ (L$^2$T$^{-1}$) signifies the hydrodynamic dispersion tensor. This tensor is made up of the following four components
\begin{linenomath*}
	\begin{equation}
	\begin{split}
	D_{xx} & = \frac{1}{\norm{ \textbf{v} }} (\alpha_\text{L} v_x^2 + \alpha_\text{T} v_y^2),\\ 
	D_{yy} & = \frac{1}{\norm{ \textbf{v} }} (\alpha_\text{L} v_y^2 + \alpha_\text{T} v_x^2),\\ 
	D_{xy} & = D_{yx} = \frac{1}{\norm{ \textbf{v} }}(\alpha_\text{L} - \alpha_\text{T}) v_{x} v_{y},
	\end{split}
	\label{dispersion}
	\end{equation}
\end{linenomath*}
where $D_{xx}$, $D_{yy}$, $D_{xy}$ and $D_{yx}$ (L$^2$T$^{-1}$) are the two principal components of the dispersion tensor and their two cross terms, respectively, $\alpha_\text{L}$ and $\alpha_\text{T}$ (L) signify the longitudinal and transverse dispersivity, respectively, $v_{x}$ and $v_{y}$ are the water flow velocities in the $x$ and $y$ direction, respectively, and $\norm{\textbf{v}} = \sqrt{v_{x}^2 + v_{y}^2}$ is the $\ell^{2}$ norm, or magnitude, of the velocity vector, $\textbf{v}$. 
\vspace*{0px}
\begin{center} [ Table \ref{table:2} NEAR HERE ] \end{center}
\vspace{0px}

To test, evaluate and benchmark the DREAM$_{\text{(KZS)}}$ algorithm, we create a reference data set as follows. We sample randomly from a multivariate uniform prior distribution of the two coordinates, $\{x_\text{s}, y_\text{s}\}$, of the contaminant source location and six release strengths, $\{s_{1},\ldots,s_{6}\}$, of the mass-loading rate step function using the ranges listed in Table \ref{table:2}. These sampled values are listed in the last column of Table \ref{table:2}, and used to construct an artificial data set by collecting, at 15 different wells within the flow domain, the simulated steady-state hydraulic heads and transient contaminant concentrations. The location of each of these wells is separately indicated in Figure \ref{fig:5} with a blue dot. We now simulate a reference run of steady-state heads and transient contaminant concentrations at each well using $\alpha_\text{L} = 0.3$ (L) and $\alpha_\text{T} = 0.03$ (L). The simulated heads and contaminant concentrations are subsequently corrupted with a zero-mean normally distributed error with standard deviations of $0.005$ (L) and $0.005$ (ML$^{-3}$), respectively. The final data set now consists of the steady-state hydraulic heads and transient contaminant concentrations at $t = \{4,5,\ldots,12\}$ (T) at the 15 wells. This equates to $n_\text{h} = 15$ measurements of the hydraulic head and $n_\text{C} = 135$ observations of the contaminant concentration. 

We now use the DREAM$_{\text{(KZS)}}$ algorithm to reconstruct the hydraulic conductivity field, point source location, and transient mass-loading rates from the measurement data. This requires a parametric definition of the hydraulic conductivity field of the rectangular flow domain. A simple Cartesian parameterization of the hydraulic conductivity field would require an excessively large number of parameters to characterize adequately the imposed spatial variability. Instead, we take advantage of the Karhunen-Lo\`{e}ve (KL) expansion and approximate the log-hydraulic conductivity field, $\mathcal{K}$, using a finite, yet relatively small, number of orthogonal basis functions \cite{zhang2004}
\begin{linenomath*}
	\begin{equation}
	\widetilde{\mathcal{K}}(\textbf{x}) = \mu_\mathcal{K}(\textbf{x}) + \sum_{i=1}^{p} \sqrt{\tau_i}s_i(\textbf{x})\beta_{i},
	\label{KL}
	\end{equation}
\end{linenomath*}
where $\mu_\mathcal{K}(\textbf{x})$ denotes the mean log-conductivity, $\textbf{x}$ represents the Cartesian coordinates of the flow domain, $s_i(\textbf{x})$ and $\tau_i, i=1,...,p,$ signify the eigenfunctions and eigenvalues of the kernel defined in equation (\ref{correlation}), and $\{\beta_{1},...,\beta_{p}\}$ are standard normal random variables, the so-called KL expansion terms. The mean square error of the reconstructed log-conductivity field, $\widetilde{\mathcal{K}}(\textbf{x})$, will go to zero in the limit of $p \to \infty$. Yet, the use of a very large number of KL terms defeats the purpose of this expansion and is not encouraged for statistical inference. Instead, one-hundred expansion terms preserve about 95\% of the variance of the ``true" log-conductivity field, $\mathcal{K}(x,y)$, i.e., $\sum_{i=1}^{p}\tau_i/\sum_{i=1}^{\infty} \tau_i \approx 0.95$. Henceforth, we characterize the conductivity field of our flow domain with $p = 100$ KL terms, $\bm{\upbeta} = \{\beta_{1},\ldots,\beta_{100}\}$.

\begin{figure}
	\centering
	\includegraphics[width=33pc]{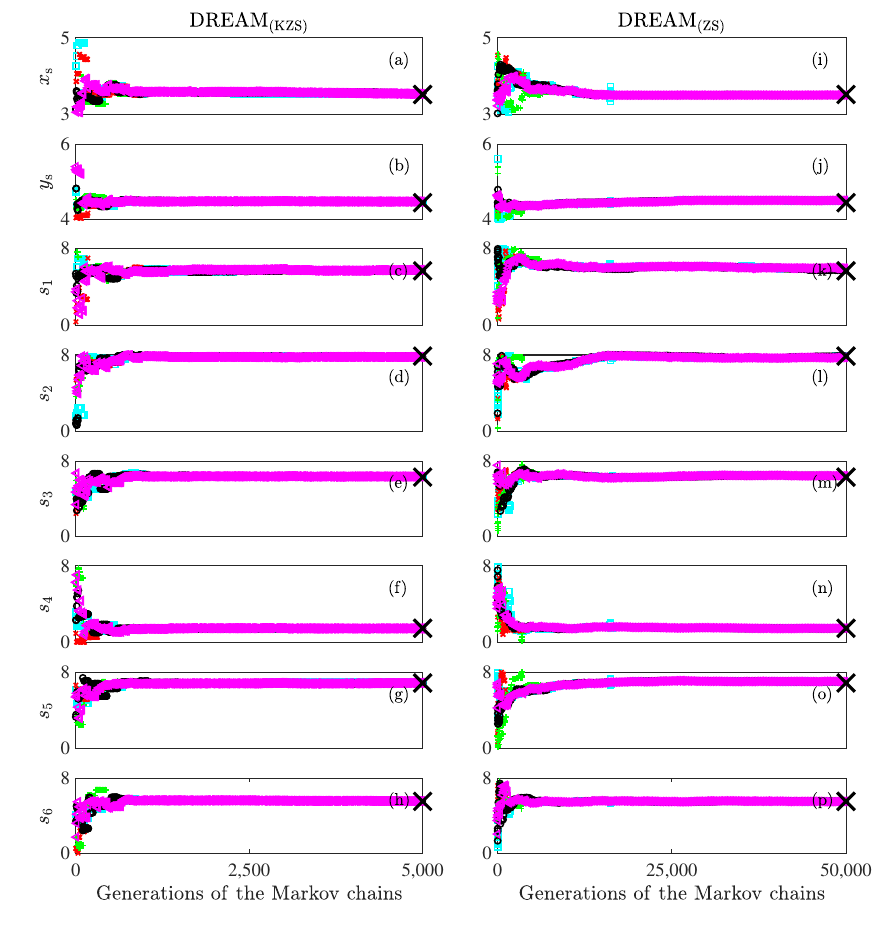}
	\caption{Trace plots of the sampled values of the point source coordinates, $\{x_\text{s}, y_\text{s}\}$, and the six coefficients, $\{s_{1},\ldots,s_{6}\}$, of the mass-loading rate step function for the DREAM$_\text{(KZS)}$ (left column) and DREAM$_\text{(ZS)}$ (right column) algorithms. Color and symbol coding are used to differentiate among the sampled chains. The ``true" parameter values are separately indicated in each panel with the black crosses.}
	\label{fig:6}
\end{figure}

Our numerical experiment now involves inference of $k = 108$ parameters, namely $p = 100$ KL terms, $\bm{\upbeta} = \{\beta_{1},\ldots,\beta_{100}\}$, two coordinates of the source location, $\{ x_\text{s}, y_\text{s} \}$, and six coefficients, $\{s_{1},\ldots,s_{6}\}$, of the mass-loading rate step function. Thus, the unknown parameters are $\bm{\uptheta} = \{\beta_{1},\ldots,\beta_{100}, x_\text{s}, y_\text{s}, s_{1},\ldots,s_{6}\}$. We use a standard Gaussian likelihood function (see equation (\ref{likelihood})) and execute the DREAM$_\text{(KZS)}$ and DREAM$_\text{(ZS)}$ algorithms with $N = 20$ chains using default values of the algorithmic variables.
\vspace*{0px}
\begin{center} [ Figure \ref{fig:6} NEAR HERE ] \end{center}
\vspace{0px}

\begin{figure}
	\centering
	\includegraphics[width=28pc]{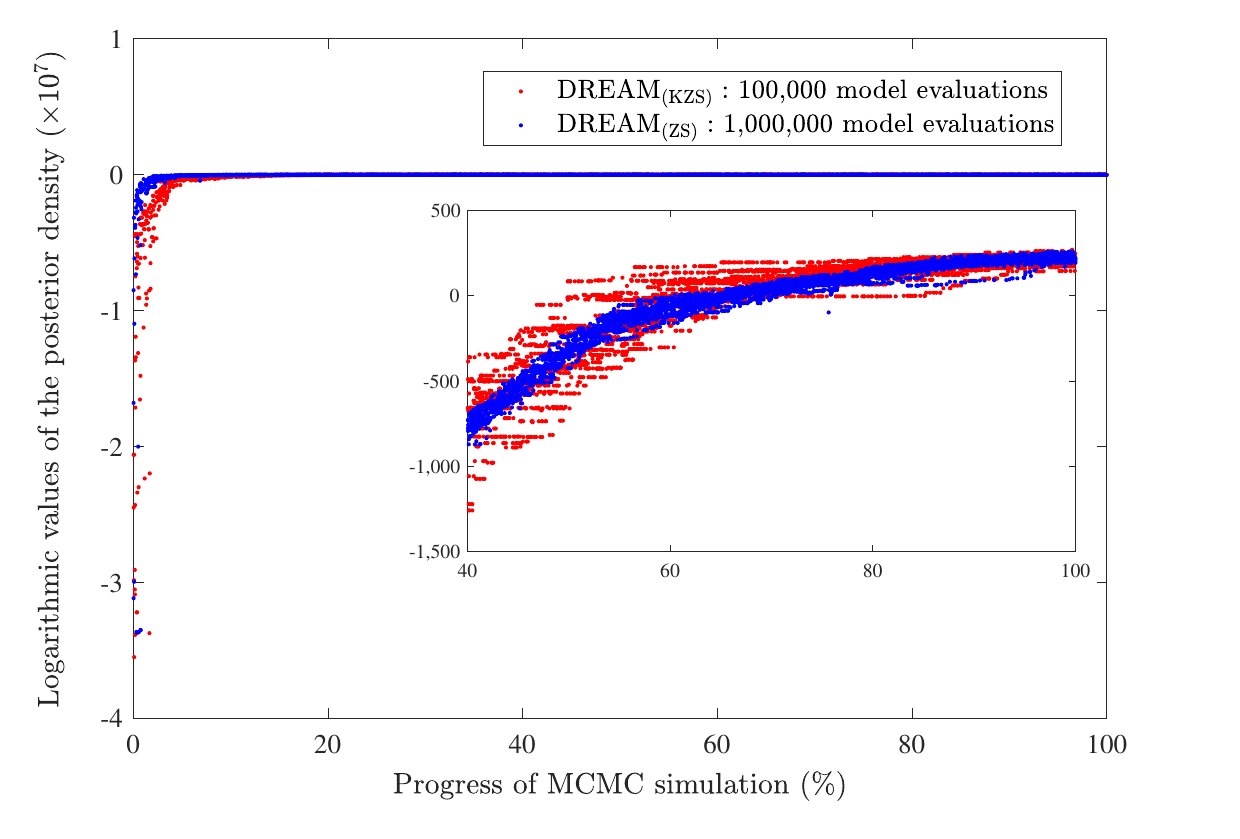}
	\caption{Evolution of the log-transformed values of unnormalized posterior density, i.e., $\mathcal{P}(\bm{\uptheta}|\widetilde{\textbf{d}})=\text{log}\bigl[p(\bm{\uptheta})L(\bm{\uptheta}|\widetilde{\textbf{d}})\bigr]$, at $\bm{\uptheta}$, sampled by both the DREAM$_\text{(KZS)}$ (red dots) and DREAM$_\text{(ZS)}$ (blue dots) algorithms.}
	\label{fig:7}
\end{figure}

\begin{figure}
	\centering
	\includegraphics[width=30pc]{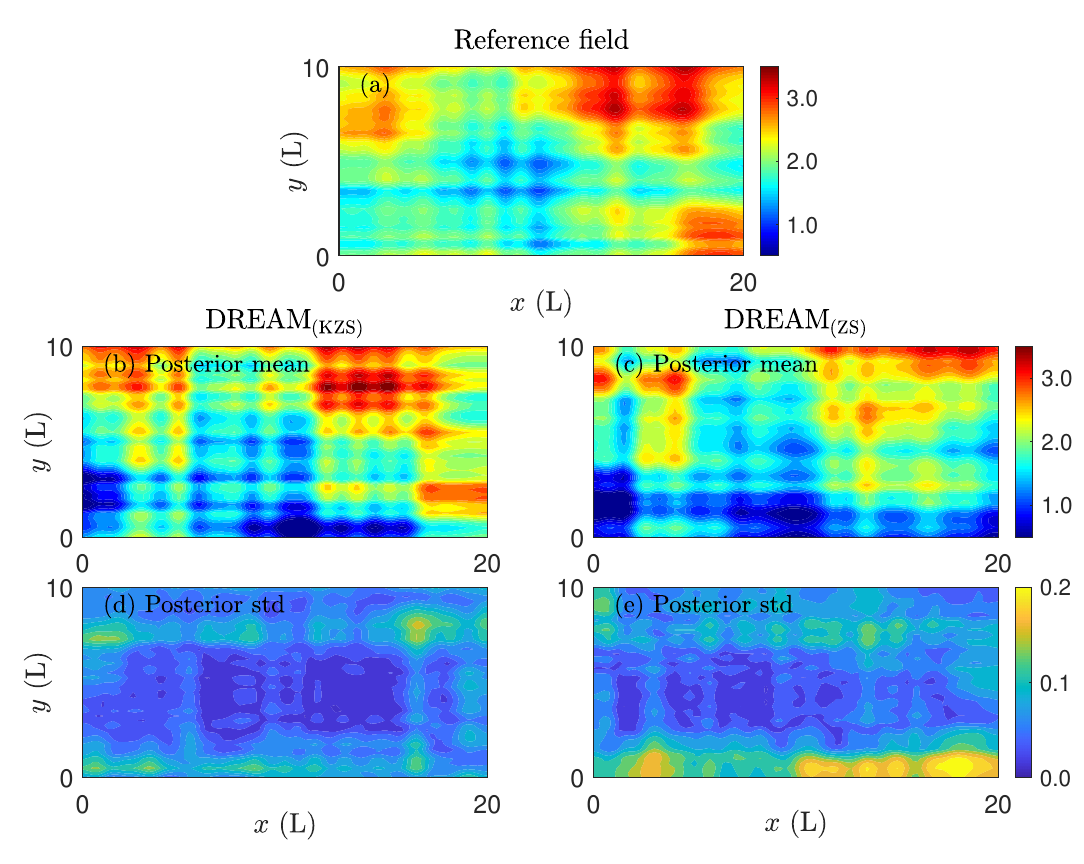}
	\caption{(a) The reference field of log-transformed conductivity, $\mathcal{K}$. (b-e) Posterior mean $\mathcal{K}$ fields (middle row) and associated standard deviations (bottom row) derived from the DREAM$_\text{(KZS)}$ (left column) and DREAM$_\text{(ZS)}$ (right column) algorithms.}
	\label{fig:8}
\end{figure}

Figure \ref{fig:6} presents trace plots of the sampled values of the point source coordinates, $\{x_\text{s}$,$y_\text{s}\}$, and the six coefficients, $\{s_{1},\ldots,s_{6}\}$, of the mass-loading rate step function for the DREAM$_\text{(KZS)}$ (left column) and DREAM$_\text{(ZS)}$ (right column) algorithms. Here we only draw five of the twenty Markov chains so that size of the image file would not be too big. The different chains are coded with different symbols and colors. The ``true" value of each parameter (i.e., the last column of Table \ref{table:2}) is separately indicated with a black cross symbol at the right-hand-side of each panel. The most important results are as follows. First, the sampled chains of both algorithms converge to the ``true" values used to generate the hydraulic head and contaminant concentration observations. Second, the location of the point source and mass-loading rate step function appear well defined with negligible posterior uncertainty compared to the width of the prior distribution. Third, and perhaps most important, the DREAM$_\text{(KZS)}$ algorithm requires far fewer model evaluations than the DREAM$_\text{(ZS)}$ algorithm to explore the target distribution.
\vspace*{0px}
\begin{center} [ Figure \ref{fig:7} NEAR HERE ] \end{center}
\vspace{0px} 

To better understand the search capabilities of the two MCMC algorithms, please consider Figure \ref{fig:7} that presents the evolution of the log-transformed values of unnormalized posterior density, i.e., $\mathcal{P}(\bm{\uptheta}|\widetilde{\textbf{d}})=\text{log}\bigl[p(\bm{\uptheta})L(\bm{\uptheta}|\widetilde{\textbf{d}})\bigr]$, at $\bm{\uptheta}$, sampled by the DREAM$_\text{(KZS)}$ (red dots) and DREAM$_\text{(ZS)}$ (blue dots) algorithms, respectively. At early stages of the search with both MCMC algorithms, the sampled chain states exhibit rather small logarithmic values (on the order of $-10^{7}$) of the posterior density. Gradually, the Markov chains move to the high probability region. During the last 60$\%$ of the MCMC simulations, $\mathcal{P}(\bm{\uptheta}|\widetilde{\textbf{d}})$ lies between -1500 and 500. Finally, both algorithms reach a similar $\mathcal{P}(\bm{\uptheta}|\widetilde{\textbf{d}})$ value of about 220. Note, that the total number of model evaluations required by the DREAM$_\text{(KZS)}$ algorithm (5,000 generations $\times$ 20 $\text{chains} = 100,000$ samples) is only one tenth of the DREAM$_\text{(ZS)}$ algorithm (50,000 generations $\times$ 20 $\text{chains} = 1,000,000$ samples). This constitutes a speed-up on the order of 10 times. With an average CPU-time of 2 seconds for each evaluation of the integrated model of MODFLOW and MT3DMS, this speed-up equates to about 500 hours reduction in CPU-cost.
\vspace*{0px}
\begin{center} [ Figure \ref{fig:8} NEAR HERE ] \end{center}
\vspace{0px}

Next, Figure \ref{fig:8} presents log-conductivity fields, $\mathcal{\widetilde{K}}(x,y)$, of the last $10,000$ posterior realizations sampled by the DREAM$_\text{(KZS)}$ (left column) and DREAM$_\text{(ZS)}$ (right column) algorithms. Figures \ref{fig:8}b and \ref{fig:8}c show the posterior mean $\mathcal{\widetilde{K}}(x,y)$ fields obtained by averaging the log-conductivity fields of the ten-thousand posterior realizations. The bottom row presents maps of the standard deviation of the ten-thousand log-conductivity fields. Compared with the reference field as shown in Figure \ref{fig:8}a, the posterior mean map of the DREAM$_\text{(KZS)}$ algorithm illuminates correctly the areas of high and low conductivity in the flow domain, yet underestimates their spatial extent. Altogether, these findings suggest that the hydraulic head and contaminant concentration measurements contain insufficient information to back out exactly the actual conductivity field. A larger and more diverse data set is warranted.

Above we have tested a nonlinear inverse problem with $k=108$ unknown model parameters. In this case, even using a state-of-the-art MCMC algorithm, i.e., DREAM$_\text{(ZS)}$, at least one million forward model evaluations are needed, which constitutes a prohibitively high computational cost. However, introducing the Kalman-inspired proposal distribution can bring a speed-up of about 10 times in simulation efficiency. The advantage of the Kalman-inspired proposal distribution is obvious in high-dimensional settings. Nevertheless, when a problem has a much larger number of unknown parameters, e.g., $k>1,000$, and this problem is highly nonlinear and complex, performance of MCMC methods will deteriorate. In this situation, one can resort to more computationally appealing methods, e.g., EnKF and its variants \cite{chen2006,crestani2013,evensen2009data}, that assume multi-Gaussian parameter and error distributions. However, these methods will not work properly when the posterior is non-Gaussian, or even has multiple modes. To address these issues, several strategies can be adopted, e.g., transforming non-Gaussian variables to be Gaussian distributed \cite{chang2010history,zhou2011approach}, or adopting a local Kalman update strategy to handle multi-modal posteriors \cite{zhang2018iterative}.

\subsection{Case Study 3: A 3-D Groundwater Model}
\begin{figure}
	\centering
	\includegraphics[width=35pc]{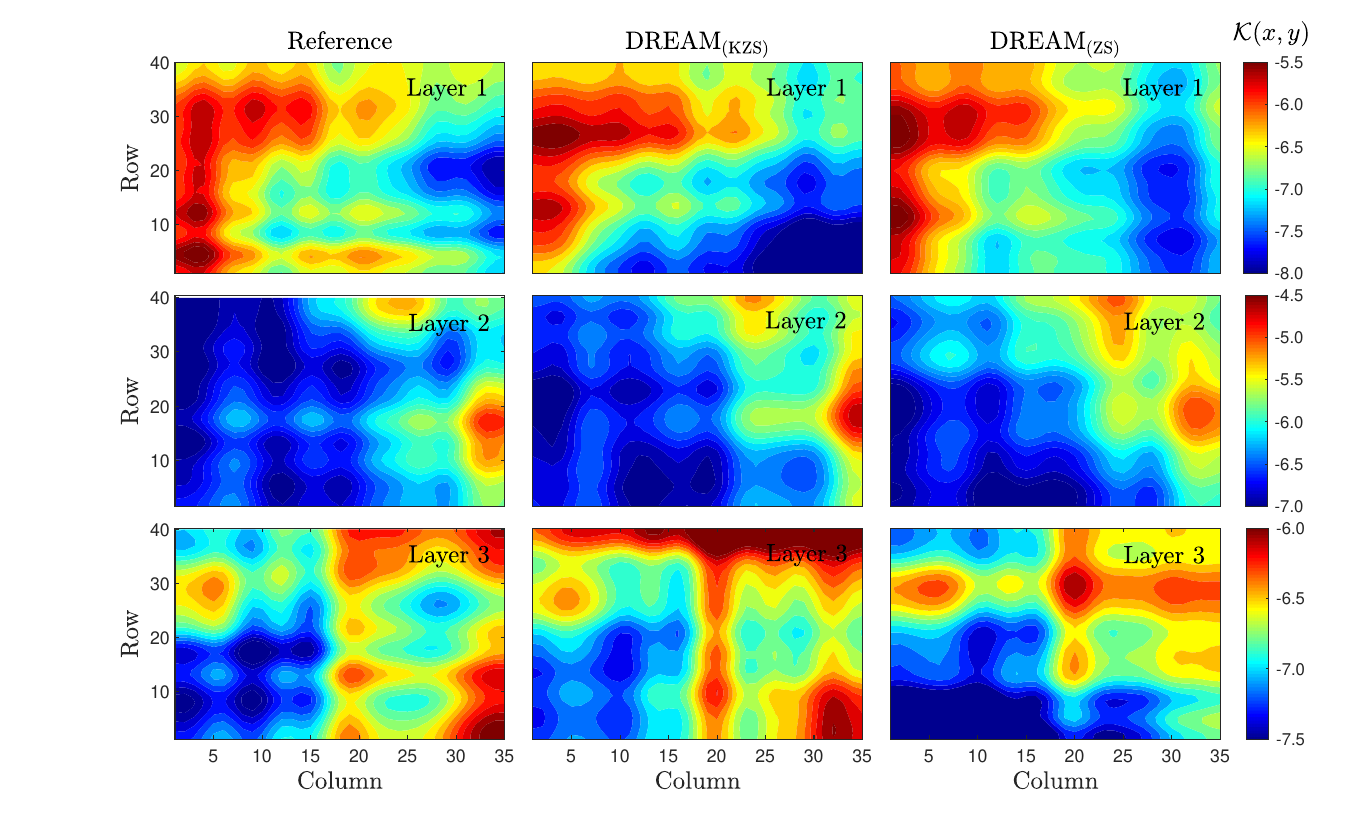}
	\caption{Reference field of the log-transformed hydraulic conductivity for each of the three layers of the aquifer (left column) and their posterior mean counterparts inferred from the DREAM$_\text{(KZS)}$ (middle column) and DREAM$_\text{(ZS)}$ (right column) algorithms using the measured stead-state hydraulic heads.}
	\label{fig:9}
\end{figure}

The third case study considers application of the DREAM$_\text{(KZS)}$ algorithm to aquifer characterization using the three-dimensional groundwater model of \citeA{fienen2013}. This model simulates three horizontal layers of 1.8, 1.4 and 1.8 m thickness that are each discretized into 35 columns and 40 rows with equidistant spacing of 2.0 and 1.5 m in the $x$ (column) and $y$ (row) direction, respectively. This equates to a domain of 70 m by 60 m by 5 m and a total of $N_\text{node} = 35 \times 40 \times 3 = 4,200$ nodes. A constant head of 60 m is prescribed at all vertical edges of the flow domain. Furthermore, an artificial well is located at row 18 and column 17, which pumps water from each horizontal layer at a constant rate of 0.01 liter/min.
\vspace*{0px}
\begin{center} [ Figure \ref{fig:9} NEAR HERE ] \end{center}
\vspace{0px}

To reduce parametric dimensionality, we use a sparse representation of the conductivity field of each layer. Here, mean and variance of the log-conductivity fields are $\mu_\mathcal{K}(\textbf{x})=-6.5$ and $\sigma^{2}_\mathcal{K}=0.5$, and correlation lengths in the $x$ and $y$ direction are $\lambda_{x} = 37.5$ m and $\lambda_{y} = 60$ m, respectively, for all the three layers. About 40 KL terms are deemed sufficient to represent the $\mathcal{K}$ field of each layer and preserve about 94\% of the field variance of each original 1,400-cell conductivity field. Hydraulic head measurements at the three layers from 81 wells located every 4 rows from row 3 to row 35, every 3 columns from column 5 to column 29, are generated from the reference fields depicted in the left column of Figure \ref{fig:9} with additive white noise, $\upepsilon \sim \mathcal{N}_n(0,0.01^2)$. We infer the $k=120$ unknown KL terms with the DREAM$_\text{(KZS)}$ and DREAM$_\text{(ZS)}$ algorithms using $N = 20$ chains and $T = 4,000$ generations. We discard the first $3,500$ generations as burn-in and use the $10,000$ samples in the last 500 generations of the joint chains to summarize the posterior estimates of the conductivity field of each discretized aquifer layer. 

\begin{figure}
	\centering
	\includegraphics[width=25pc]{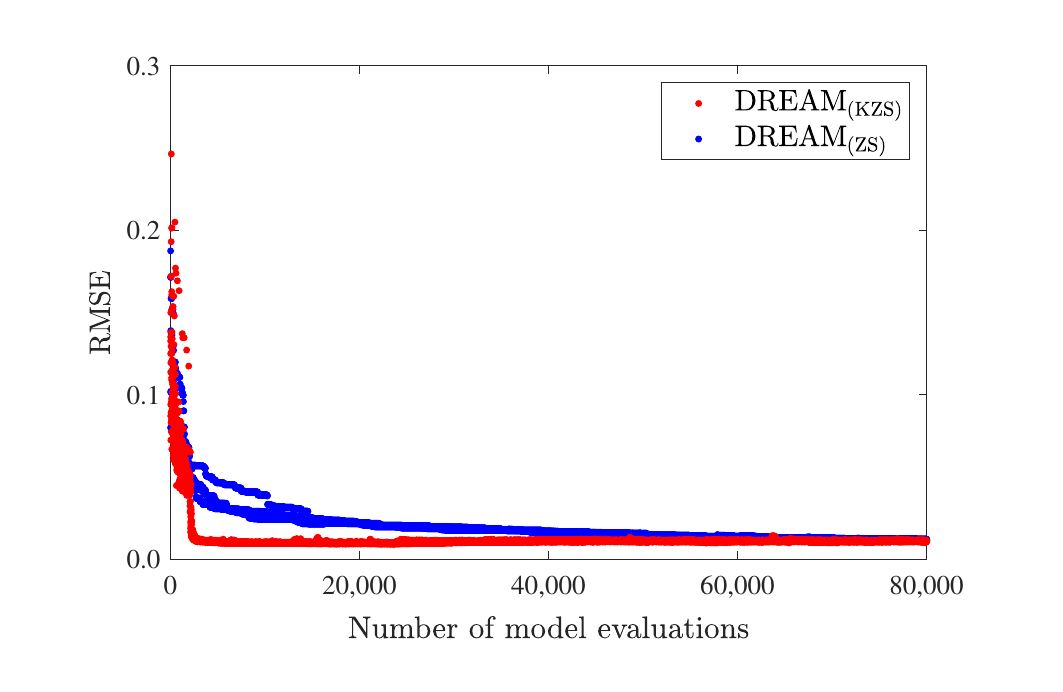}
	\caption{Evolution of the sampled RMSE values between the MODFLOW simulated and observed steady-state hydraulic heads. All chains are plotted. Color coding in red and blue differentiates between the chains sampled by the DREAM$_\text{(KZS)}$ and DREAM$_\text{(ZS)}$ algorithms, respectively.}
	\label{fig:10}
\end{figure}

Figure \ref{fig:9} presents the posterior mean log-conductivity fields derived from the DREAM$_\text{(KZS)}$ (middle column) and DREAM$_\text{(ZS)}$ (right column) algorithms. The mean log-conductivity fields derived from both algorithms capture quite well the main patterns of the reference fields of the three layers (left column). Some discrepancies are visible, but appear relatively minor. These results are encouraging but do not convey anything about the efficiency of the two MCMC algorithms. We therefore proceed with analysis of the convergence properties of the DREAM$_\text{(KZS)}$ and DREAM$_\text{(ZS)}$ algorithms.
\vspace*{0px}
\begin{center} [ Figure \ref{fig:10} NEAR HERE ] \end{center}
\vspace{0px}

In Figure \ref{fig:10}, we present trace plots of the root-mean-square error (RMSE) in the $N = 20$ Markov chains between the MODFLOW simulated and measured steady-state heads at the 81 wells. We use color coding in red and blue for the DREAM$_\text{(KZS)}$ and DREAM$_\text{(ZS)}$ algorithms, respectively. The results in this figure confirm our earlier conclusions. The DREAM$_\text{(KZS)}$ algorithm converges at a much faster pace to the target distribution. Indeed, whereas the DREAM$_\text{(ZS)}$ algorithm requires about 60,000 - 80,000 function evaluations to sample values of the RMSE on the order of the Gaussian measurement error of 0.01, the DREAM$_\text{(KZS)}$ algorithm needs only about 2,500 - 3,000 MODFLOW evaluations to minimize the RMSE. This equates to a speed-up on the order of 20 - 30 times. This gives the DREAM$_\text{(KZS)}$ algorithm sufficient opportunity to sample thoroughly the target distribution and summarize adequately the posterior moments of the quasi three-dimensional conductivity distribution of the aquifer. 

In this case study, although the number of unknown model parameters ($k=120$) is slightly larger than that in the second case study ($k=108$), much fewer model evaluations are needed in the MCMC simulations. Using the DREAM$_\text{(ZS)}$ algorithm, to obtain acceptable results, the total numbers of model evaluations in the second and present case studies are at least 1,000,000 and 60,000, respectively. Using the DREAM$_\text{(KZS)}$ algorithm, the corresponding numbers are about 100,000 and 2,500. This comparison indicates that the convergence speed of an MCMC algorithm is determined not only by the number of model parameters, but also by the nonlinearity and complexity of the problem. In many situations, the high-dimensionality of unknown parameters (e.g., hydraulic conductivity field) comes from the discretization of the underlying function. Using a standard MCMC algorithm, the convergence speed degrades significantly with mesh refinement. This kind of MCMC algorithm is thus called dimension-dependent. To address this kind of dimensionality issue, dimension-independent MCMC methods have been developed, e.g., the preconditioned Crank-Nicolson algorithm proposed by \citeA{cotter2013mcmc} that preserves the reference measure with a Crank-Nicolson discretization of the underlying function, and the dimension-independent and likelihood-informed algorithm proposed by \citeA{cui2016dimension} that constructs a global likelihood-informed subspace to capture essential features of the posterior distribution. The dimension-independent MCMC methods require the distributions exhibit certain special structure. Yet, sampling generic high-dimensional target distributions with MCMC methods is still a challenging problem.

\subsection{Case Study 4: Inverse Problems with Multi-modal Posteriors}
	
To demonstrate the performance of the DREAM$_\text{(KZS)}$ algorithm in solving nonlinear and ill-posed inverse problems, three examples with multi-modal posterior distributions are further tested below.

\begin{figure}
	\centering
	\includegraphics[width=30pc]{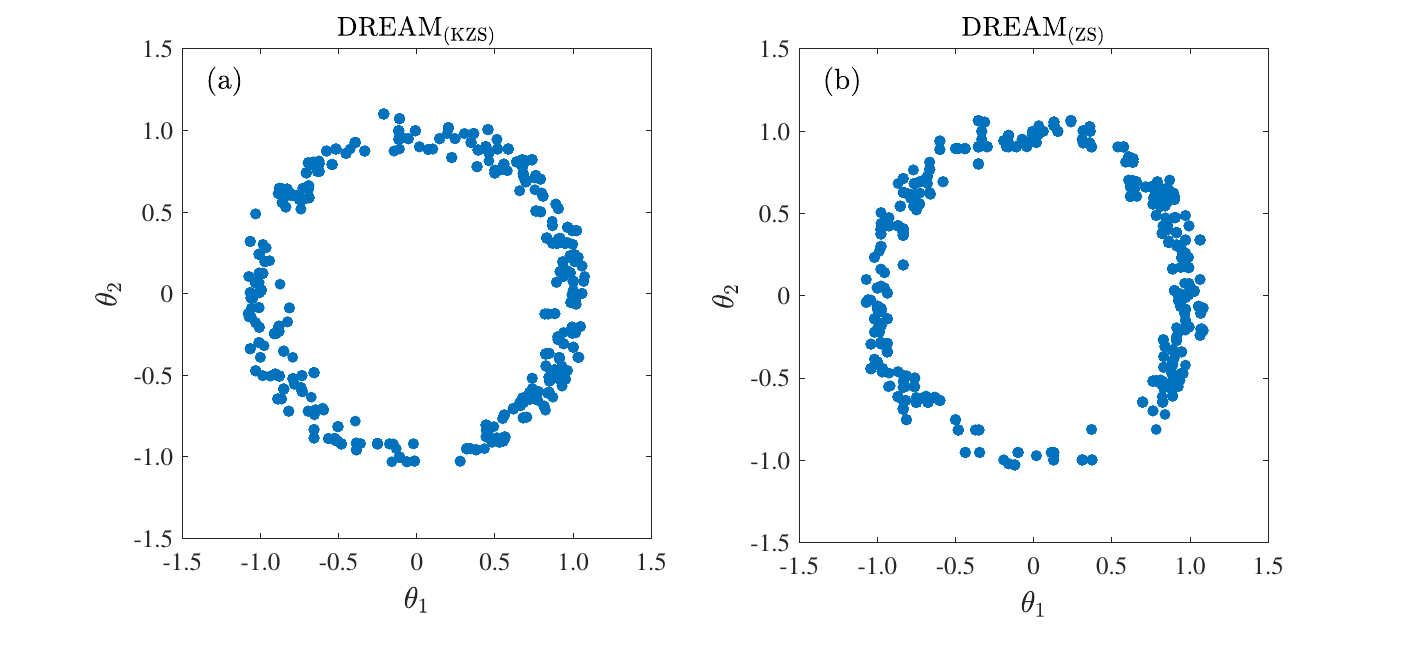}
	\caption{Bivariate scatter plots of posterior parameter samples obtained by (a) DREAM$_\text{(KZS)}$ and (b) DREAM$_\text{(ZS)}$.}
	\label{fig:11}
\end{figure}

The first example is very simple, yet highly ill-posed. The underlying function has the following form
\begin{linenomath*}
	\begin{equation}
	d=\theta_{1}^{2}+\theta_{2}^{2},
	\label{simple_mul_1}
	\end{equation}
\end{linenomath*}
where $\theta_{1} \sim \mathcal{U}(-2,2)$ and $\theta_{2} \sim \mathcal{U}(-2,2)$ are all uniformly distributed in the prior. The measurement is $\widetilde{d}=1$ with measurement error, $\epsilon \sim \mathcal{N}\left(0,0.1^{2}\right)$. Here we are trying to infer the joint posterior distribution of two parameters, $\bm{\uptheta}=\left\{\theta_{1}, \theta_{2}\right\}$, from a scalar measurement, $\widetilde{d}$. It is evident that there exist an infinite number of parameter combinations that can fit the measurement well. Figure \ref{fig:11} shows the bivariate scatter plots of posterior parameter samples obtained by both the DREAM$_\text{(KZS)}$ and DREAM$_\text{(ZS)}$ algorithms. Here both algorithms evolve $N=3$ Markov chains with $T=1,000$ samples in each chain. The last 300 samples in each of the three chains are used to draw Figure \ref{fig:11}. It is found that both algorithms obtain reasonable results.
\vspace*{0px}
\begin{center} [ Figure \ref{fig:11} NEAR HERE ] \end{center}
\vspace{0px}

\begin{figure}
	\centering
	\includegraphics[width=30pc]{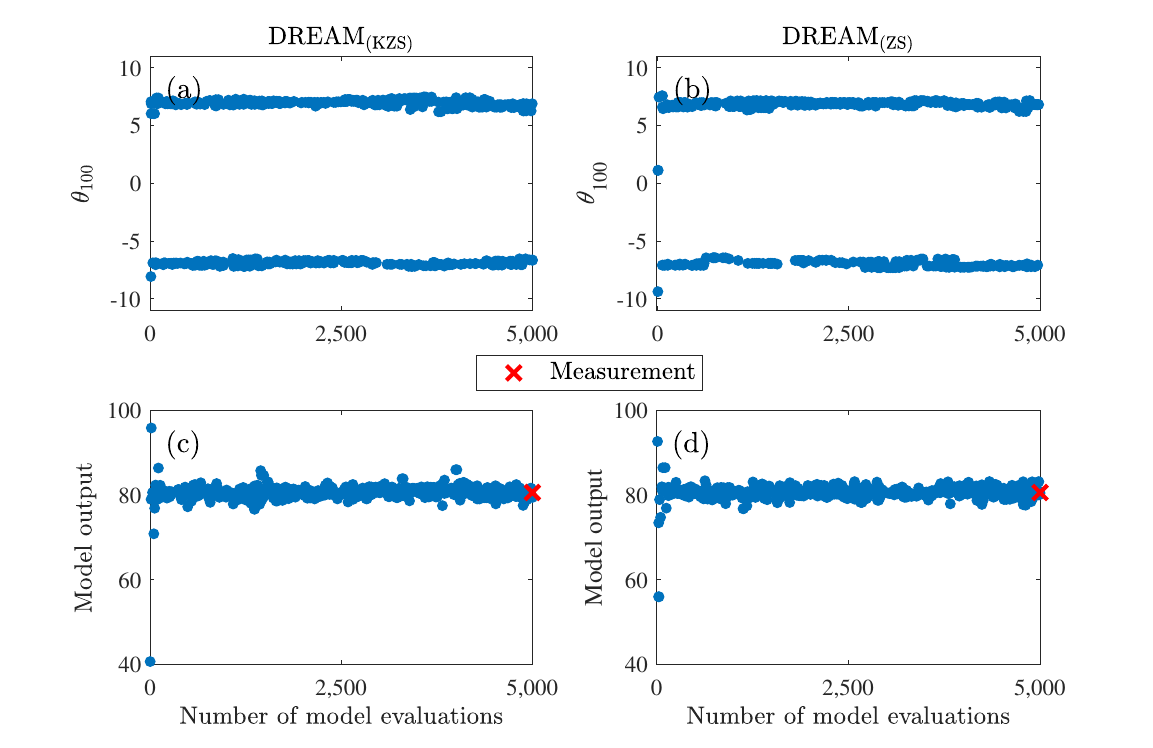}
	\caption{Trace plots of the sampled values of (a-b) $\theta_{100}$ and (c-d) the model output. The left column is for the DREAM$_\text{(KZS)}$ algorithm, and the right column is for the DREAM$_\text{(ZS)}$ algorithm. The measurement $\widetilde{d}$ is represented by a red cross in each panel at the bottom row.}
	\label{fig:12}
\end{figure}

Then we extend the above example to a high-dimensional setting, which turns the underlying function to
\begin{linenomath*}
	\begin{equation}
	d=\theta_{1}^{2}+\theta_{2}^{2}+\cdots+\theta_{100}^{2},
	\label{simple_mul_2}
	\end{equation}
\end{linenomath*}
where $\left\{\theta_{1}, \dots, \theta_{99}\right\}$ are all uniformly distributed with $\mathcal{U}(0,2)$, and $\theta_{100} \sim \mathcal{U}(-10,10)$. There is still only a scalar measurement, whose value is $\widetilde{d}=80.62$, and the error is normally distributed, $\epsilon \sim \mathcal{N}\left(0,1^{2}\right)$. Compared to the previous example, inferring the joint distribution of 100-dimensional parameters, $\bm{\uptheta}=\left\{\theta_{1},\cdots,\theta_{100}\right\}$, from a scalar measurement, $\widetilde{d}$, becomes more challenging. As $\theta_{100}$ has a wide prior range, the marginal posterior distribution of $\theta_{100}$ can be obviously bi-modal, as both $\theta_{100}$ and $-\theta_{100}$ have the same effect on the function output. As shown in Figure \ref{fig:12}, both the DREAM$_\text{(KZS)}$ and DREAM$_\text{(ZS)}$ algorithms (here $N=5$ and $T=1,000$) can identify the bimodality of $\theta_{100}$ and fit the measurement quite well.
\vspace*{0px}
\begin{center} [ Figure \ref{fig:12} NEAR HERE ] \end{center}
\vspace{0px}
\vspace*{0px}
\begin{center} [ Table \ref{table:3} NEAR HERE ] \end{center}
\vspace{0px}

\begin{figure}
	\centering
	\includegraphics[width=40pc]{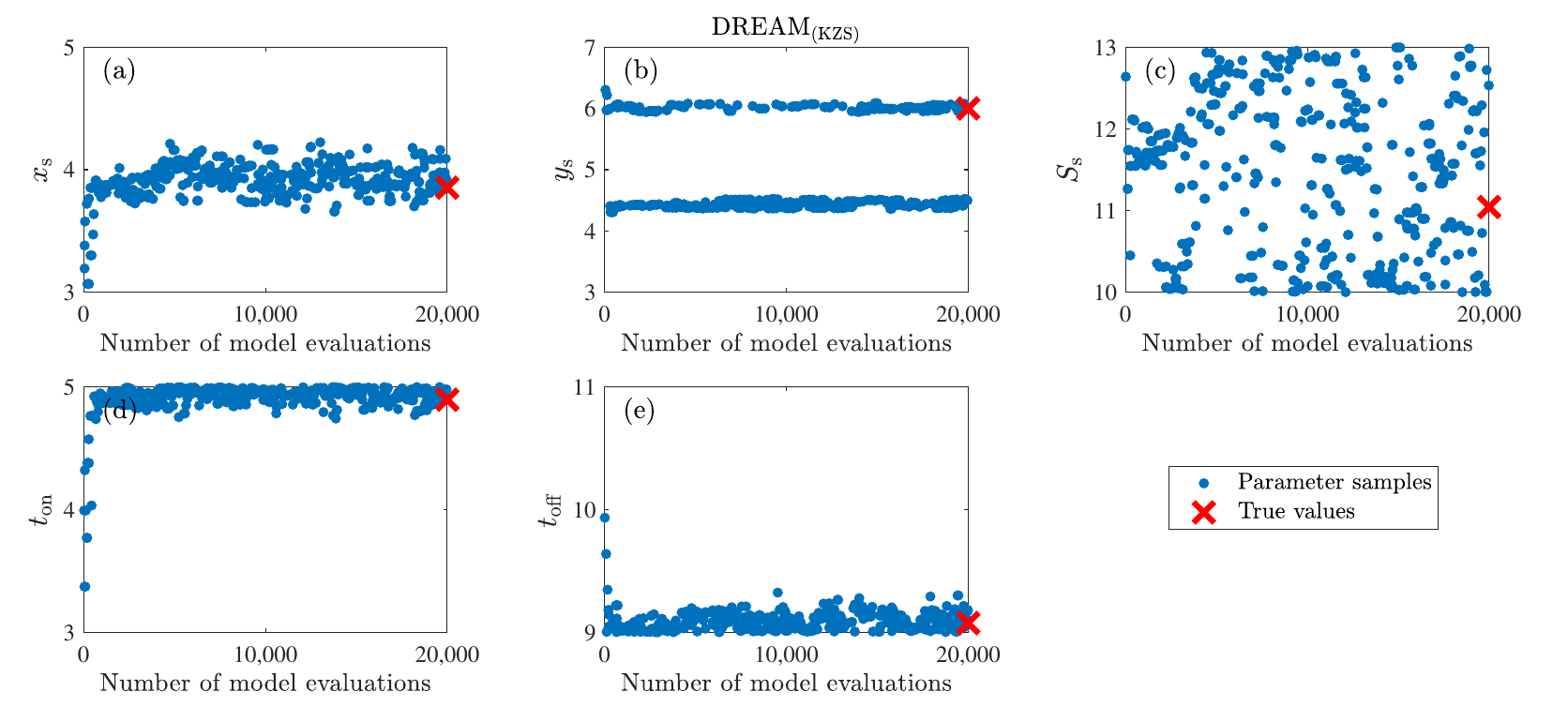}
	\caption{Trace plots of the sampled values of the five contaminant source parameters obtained by the DREAM$_\text{(KZS)}$ algorithm. The ``true” parameter values are separately indicated in each panel with the red crosses.}
	\label{fig:13}
\end{figure}

Finally, we test a groundwater contaminant source identification problem with multiple modes in the posterior. The model settings are the same as the example tested in section 3.2, except that here we consider a homogeneous conductivity field with a known value of $K=8\left(\mathrm{LT}^{-1}\right)$, and a different parameterization of the contaminant source. The source located at $\left(x_{\mathrm{s}}, y_{\mathrm{s}}\right)$ (L) starts to release from $t=t_{\mathrm{on}}$ (T) with a constant mass-loading rate of $S_{\mathrm{s}}$ (MT$^{-1}$) until $t=t_{\mathrm{off}}$ (T). Thus, there are five unknown parameters, i.e., $\bm{\uptheta} =\left\{x_{\mathrm{s}}, y_{\mathrm{s}}, S_{\mathrm{s}}, t_{\mathrm{on}}, t_{\mathrm{off}}\right\}$, whose prior ranges (uniform distributions) and “true” values are provided in Table \ref{table:3}. Concentration measurements are collected at a single well (the blue dot in Figure S9) at $t=\{6,8,10,12,14\} (\mathrm{T})$. The measurement error is normally distributed, $\epsilon \sim \mathcal{N}\left(0,0.01^{2}\right)$. Then we implement the DREAM$_\text{(KZS)}$ and DREAM$_\text{(ZS)}$ algorithms to infer the posterior distribution of the five unknown parameters, respectively. Here, for both algorithms, there are $N=4$ parallel chains, and each chain has $5,000$ samples. Trace plots of the sampled values of the five contaminant source parameters, i.e., $\bm{\uptheta} =\left\{x_{\mathrm{s}}, y_{\mathrm{s}}, S_{\mathrm{s}}, t_{\mathrm{on}}, t_{\mathrm{off}}\right\}$, obtained by the DREAM$_\text{(KZS)}$ and DREAM$_\text{(ZS)}$ algorithms are depicted in Figures \ref{fig:13} and S10, respectively. Both algorithms can identify the bimodal posterior distribution of $y_{\mathrm{s}}$. From the multivariate $\widehat{R}^{k}$-diagnostic that monitors the convergence of the Markov chains (Figure S11), it is found that the DREAM$_\text{(KZS)}$ algorithm converges to its stationary regime slightly faster than the DREAM$_\text{(ZS)}$ algorithm.
\vspace*{0px}
\begin{center} [ Figure \ref{fig:13} NEAR HERE ] \end{center}
\vspace{0px}

The above three examples demonstrate that the proposed method can still work properly in multi-modal cases. When dealing with a nonlinear, non-Gaussian inverse problem, the simulation results of EnKF and it variants will deteriorate significantly, as the updated ensemble of states cannot approximate the complex posterior distribution accurately. However, the Kalman-inspired proposal used in the MCMC simulation can still have a positive effect. In the Kalman-inspired proposal, we only generate a candidate at a time. Although the jump from the current state $\bm{\uptheta}_{(t-1)}$ to the candidate $\bm{\uptheta}_\text{p}$ may be not optimal, it can still provide some information about the high posterior density region, and thus shorten burn-in. Even in the worst case that most of the Kalman-inspired candidates are rejected (very unlikely though), as we only use the Kalman-inspired proposal during a relatively short burn-in period (e.g., the first 30 \% of the MCMC simulation) with a small selection probability (e.g., $p_\text{K}=0.3$), only less than 9\% of the total computational cost will be wasted because of introducing the Kalman-inspired proposal distribution. The complex posterior can still be explored with other proposal distributions (e.g., the parallel direction and snooker proposal distributions). If we monitor the rejection rate of the Kalman-inspired proposal distribution and find a very high value, we can actually stop using the Kalman-inspired proposal early to avoid further computational waste. Nevertheless, in complicated, multi-modal cases, the proposed method still faces the risk of missing secondary modes. To improve the capacity of the Kalman-inspired proposal in multi-modal cases, we can adopt the modified Kalman-inspired proposal suggested in the last paragraph of section 2.3.1. The modified method borrows ideas from a new ensemble-smoother-based method \cite{zhang2018iterative} that is specifically designed to solve inverse problems with (possible) multi-modal parameter distributions. 

\section{Discussion and Conclusions}
MCMC methods have found widespread application and use to approximate the posterior distribution. Such methods generate a random walk through the parameter space and successively visit solutions with frequency proportional to the density of the underlying target distribution. The speed with which MCMC methods converge to the stationary distribution, however, deteriorates rapidly with increasing target dimensionality.

The power and usefulness of the Kalman analysis step has been demonstrated time and again with application to state estimation in real-time forecasting studies. The use of the analysis state enhances considerably the short-term predictive skill of computer simulation models. The analysis step in data assimilation methods can facilitate parameter estimation as well, which makes possible the application of data assimilation methods such as EnKF and its variants to solving high-dimensional inverse problems. The so-obtained posterior parameter distribution can, at best, only roughly approximate the target distribution.

This paper introduces a Kalman-inspired proposal distribution to improve the efficiency of posterior exploration using MCMC methods. This new proposal distribution exploits the cross-covariances of model parameters, measurements and model outputs, and generates candidate states much alike the analysis step in the Kalman filter. The Kalman-inspired proposal distribution was embedded in the DREAM algorithm, and this new sibling of the DREAM family of MCMC methods coined the DREAM$_\text{(KZS)}$ algorithm. As the Kalman-inspired proposal distribution is asymmetric, its use is restricted to a relatively short burn-in period, after which a mix of parallel direction and snooker candidate states are used to evolve the chains in the DREAM$_\text{(KZS)}$ algorithm. Diminishing adaptation guarantees that the sampled chains converge to the exact target distribution. Numerical experiments with watershed and aquifer models confirm that the Kalman-inspired proposal distribution enhances considerably the efficiency of posterior exploration. Specifically, we observe a speed-up on the order of 20 - 30 times for a three-dimensional groundwater model with more than one-hundred unknown parameters. Although we combine the Kalman-inspired idea with the DREAM algorithm, the new proposal distribution is not tied to a specific MCMC algorithm, but can be conveniently embedded in any adequate MCMC method.

When the inverse problem is ill-posed and the posterior is multi-modal, introducing the Kalman-inspired proposal distribution into an MCMC algorithm can still help. Furthermore, there have been some strategies developed for EnKF and its variants to solve non-Gaussian inverse problems, e.g., transforming non-Gaussian variables to be Gaussian distributed \cite{chang2010history,zhou2011approach} or adopting a local Kalman update strategy to handle multi-modal posteriors \cite{zhang2018iterative}. These strategies can also be used in the Kalman-inspired proposal to gain further strength. On a more theoretical note, it may be desirable to enforce symmetry of the Kalman-inspired proposal distribution. This would deteriorate at least somewhat sampling efficiency but make possible its application to the entire chain generations simulated by the MCMC algorithm. These issues will be addressed in our future work.

\acknowledgments
Computer codes and data used are available at\\
$\texttt{https://www.researchgate.net/publication/318645382\_MATLAB\_codes\_of\_DREAM\_KZS}$.\\
This work is supported by the National Key Research
and Development Program of China (grant 2018YFC1800303), National Natural Science Foundation of China (grants 41807006 and 41771254) and China Postdoctoral Science Foundation funded project (grant 2018M630680). The authors would like to thank the Editor and anonymous reviewers for their constructive comments and suggestions, which significantly improve the quality of this work. The authors would also like to thank Marko Laine from Finnish Meteorological Institute for providing the MATLAB codes of the Adaptive Metropolis and Delayed Rejection Adaptive Metropolis algorithms.


%
%

\bibliography{myref}

%
%
%
%
%

\clearpage

\begin{table}
	\caption{Description of the hmodel parameters and their minimum and maximum values that define the prior uncertainty ranges. The last column lists the ``true" parameter values used in our numerical experiment.}
	\centering
	\begin{tabular}{c c c c c c}
		\hline
		Parameter & Symbol & Min. & Max. & Unit & True \\
		\hline
		Maximum interception            &   $I_\text{max}$       & 0.5   & 10    & mm    & 3.84    \\
		Soil water storage capacity     &   $S_\text{max}$       & 10    & 1000  & mm    & 776.40  \\
		Maximum percolation rate        &   $Q_\text{max}$       & 0     & 100   & mm/d  & 26.60   \\
		Evaporation parameter           &   $\alpha_\text{E}$    & 0     & 100   & -     & 61.61   \\
		Runoff parameter                &   $\alpha_\text{F}$    & -10   & 10    & -     &-4.18    \\
		Time constant, fast reservoir   &   $K_\text{F}$         & 0     & 10    & days  & 6.01    \\
		Time constant, slow reservoir   &   $K_\text{S}$         & 0     & 150   & days  & 111.67  \\
		\hline
		\label{table:1}
	\end{tabular}
\end{table}

\begin{table}
	\caption{Description of the eight contaminant source parameters of our contaminant transport model, and the ranges of their multivariate uniform prior distribution. The last column lists the ``true" values of the parameters used in our numerical experiment.}
	\centering
		\begin{tabular}{c c c c c c}
			\hline
			Parameter & Symbol & Min. & Max. & Unit & True \\
			\hline
			$x$-coordinate of spill                 &   $x_\text{s}$   & 3     & 5    & L          & 3.52   \\
			$y$-coordinate of spill                 &   $y_\text{s}$   & 4     & 6    & L          & 4.44   \\
			Release strength of the first segment   &   $s_1$          & 0     & 8    & MT$^{-1}$  & 5.69   \\
			Release strength of the second segment  &   $s_2$          & 0     & 8    & MT$^{-1}$  & 7.88   \\
			Release strength of the third segment   &   $s_3$          & 0     & 8    & MT$^{-1}$  & 6.31   \\
			Release strength of the fourth segment  &   $s_4$          & 0     & 8    & MT$^{-1}$  & 1.49   \\
			Release strength of the fifth segment   &   $s_5$          & 0     & 8    & MT$^{-1}$  & 6.87   \\
			Release strength of the sixth segment   &   $s_6$          & 0     & 8    & MT$^{-1}$  & 5.55   \\
			\hline
			\label{table:2}
		\end{tabular}
\end{table}

\begin{table}
	\caption{Description of the five contaminant source parameters in the multi-modal case, and the ranges of their multivariate uniform prior distribution. The last column lists the ``true" values of the parameters used in our numerical experiment.}
	\centering
	\begin{tabular}{c c c c c c}
		\hline
		Parameter & Symbol & Min. & Max. & Unit & True \\
		\hline
		$x$-coordinate of spill                 &   $x_\text{s}$   & 3     & 5    & L          & 3.85   \\
		$y$-coordinate of spill                 &   $y_\text{s}$   & 3     & 7    & L          & 6.00   \\
		Constant release strength               &   $S_\text{s}$   & 10    & 13   & MT$^{-1}$  & 11.04  \\
		Start time of contaminant release       &   $t_\text{on}$  & 3     & 5    & T          & 4.90   \\
		End time of contaminant release         &   $t_\text{off}$ & 9     & 11   & T          & 9.08   \\
		\hline
		\label{table:3}
	\end{tabular}
\end{table}

\clearpage

\newpage

\newpage

\newpage

\newpage

\newpage

\newpage

\end{document}